\documentclass[12pt,reqno]{amsart}
\usepackage{etex}
\usepackage{amsmath,amsfonts,amssymb}
\usepackage{mathrsfs}
\usepackage{array}
\usepackage{amsthm,youngtab}
\usepackage{graphicx}
\usepackage{caption}
\usepackage{resizegather}
\usepackage{enumerate,enumitem}
\usepackage[textwidth=16cm, hmarginratio=1:1]{geometry}
\usepackage{adjustbox}
\usepackage{tikz-cd}
\usepackage{rotating} 
\usepackage{rotating}
\usepackage{mathdots}
\usepackage[all]{xy}
\usepackage{xcolor}
\usepackage[titletoc,title]{appendix}
\usepackage{amscd}
\usepackage[colorlinks,backref=page]{hyperref} 
\usepackage{float}
\usepackage{makecell}
\usepackage{listings}
\usepackage{tcolorbox}
\allowdisplaybreaks[4]

\hypersetup{
bookmarksnumbered,
pdfstartview={FitH},
breaklinks=true,
linkcolor=blue,
urlcolor=blue,
citecolor=blue,
bookmarksdepth=2
}

\theoremstyle{plain}
   \newtheorem{theorem}{Theorem}[section]
   \newtheorem{proposition}[theorem]{Proposition}
   \newtheorem{lemma}[theorem]{Lemma}
   \newtheorem{corollary}[theorem]{Corollary}
   \newtheorem{conjecture}[theorem]{Conjecture}
\theoremstyle{definition}
   \newtheorem{definition}[theorem]{Definition}
   \newtheorem{example}[theorem]{Example}
   \newtheorem{problem}[theorem]{Problem}
   
   \newtheorem{remark}[theorem]{Remark}

\numberwithin{equation}{section}

\newcommand{\g}{\mathfrak{g}}
\newcommand{\cc}{\textbf{coe}}
\newcommand{\hw}{\textbf{hw}}

\begin{document}
\begin{sloppypar}

\title[Classification of real modules]{Classification of real modules in monoidal categorifications of cluster algebras}
\author{Bing Duan}
\address[Bing Duan]{School of Mathematics and Statistics, Lanzhou University, Lanzhou, 730000 P. R. China}
\email{duanbing@lzu.edu.cn}
\author{Ralf Schiffler}
\address[Ralf Schiffler]{Department of Mathematics, University of Connecticut, Storrs, CT 06269--1009, USA}
\email{schiffler@math.uconn.edu}
\thanks{The first author was supported by Lanzhou University (No. 561120205). The second author was supported by the NSF grant DMS-2348909.}

\maketitle

\begin{abstract}
In this paper, we propose a conjectural formula for the highest $\ell$-weight monomial of an arbitrary real module over a simply-laced quantum affine algebra. We verify the conjecture under a multiplicative reachability condition, answering the Hernandez--Leclerc classification problem in monoidal categorifications of cluster algebras under this condition. Moreover, we introduce the notion of cluster modules, generalizing Kirillov--Reshetikhin modules and Hernandez--Leclerc modules as special cases. We prove that cluster modules are reachable real modules, and obtain a system of equations governing $q$-characters of the prime cluster modules, providing a natural generalization of both the classical T-system relations for Kirillov--Reshetikhin modules and the exchange relations for Hernandez--Leclerc modules.

\hspace{0.15cm}


\hspace{0.15cm}

\end{abstract}

\addtocontents{toc}{\protect\setcounter{tocdepth}{1}}
\tableofcontents

\section{Introduction}

Cluster algebras were introduced by Fomin and Zelevinsky \cite{FZ02} to provide a combinatorial framework for studying canonical bases and their duals in quantum groups, as well as total positivity in algebraic groups. Subsequently, these algebras have found diverse applications across mathematics and physics.

A powerful approach to studying cluster algebras is through \textit{categorification}. This strategy replaces core components of cluster algebras by analogous categorical structures with similar properties, we refer to \cite{Rei10} for a comprehensive exposition. Two significant categorifications exist, that is, additive categorification and monoidal categorification.

Certain cluster algebras admit both additive and monoidal categorifications simultaneously. Cluster categories for acyclic quivers (coefficient-free) were initially developed in \cite{BMRRT06}, with the type $A_n$ case treated separately in \cite{CCS06}. For any acyclic quiver $Q$, let $\mathcal{C}_Q$ be the cluster category and $\mathcal{A}(Q)$ the cluster algebra. A foundational result states that there exists a canonical bijection between the isomorphism classes of indecomposable rigid objects in $\mathcal{C}_Q$ and the cluster variables in $\mathcal{A}(Q)$. It was generalized beyond acyclic quivers in \cite{Ami09}. In \cite{DS25}, we established a direct link between the two categorifications for the cluster algebra associated to Hernandez--Leclerc's level-one subcategory. Recent advances in this interplay have attracted growing attention \cite{BFL24,Con24,F24,Qin20}.

Let $\mathfrak{g}$ be a simply-laced simple Lie algebra over $\mathbb{C}$ with Dynkin index set $I$, and $U_q(\widehat{\mathfrak{g}})$ its associated quantum affine algebra. The complete classification of simple $U_q(\widehat{\mathfrak{g}})$-modules was established in \cite{CP91,CP94} via Drinfeld polynomial parametrization. Combined with $q$-character theory \cite{FR98}, it follows that a simple module $L(m)$ is uniquely parameterized by its highest $\ell$-weight monomial $m$. Every simple module arises as the head of an ordered tensor product of fundamental modules. Moreover, geometric frameworks for studying these modules have been systematically developed through seminal work in \cite{Nak01,VV02}.

Let $\mathscr{C}_{\mathfrak{g}}$ denote the category of type 1 finite-dimensional $U_q(\widehat{\mathfrak{g}})$-modules. The Hopf algebra structure of $U_q(\widehat{\mathfrak{g}})$ endows $\mathscr{C}_{\mathfrak{g}}$ with a natural monoidal structure. For any height function $\xi: I \to \mathbb{Z}$, determining a Dynkin quiver $Q_\xi$, and integer $l \geq 1$, there exists a monoidal subcategory $\mathscr{C}^{\leq \xi}_l$ whose Grothendieck ring $K_0(\mathscr{C}^{\leq \xi}_l)$ admits a cluster algebra structure with initial seed $(\mathbf{z}^{\leq \xi}_l, \Gamma^{\leq \xi}_l)$, see \cite{HL10,HL16,Bit21,DS25,FHOO23}.

A simple $U_q(\widehat{\mathfrak{g}})$-module is \textit{real} if its tensor square remains simple, and \textit{imaginary} otherwise \cite{Lec03}. A simple module is \textit{prime} if it cannot be expressed as a tensor product of non-trivial modules. The study of real and imaginary modules is intimately connected to the dual canonical basis \cite{Lec03}. In \cite[Section 5.2.6]{HL21}, Hernandez and Leclerc proposed the following open problem:

\begin{problem}\label{Hernandez and Leclerc's classification open question}
Classify real modules in terms of their highest $\ell$-weight monomials.
\end{problem}

\hspace{-0.1cm}

Lapid and Minguez \cite{LM18} classified all real modules satisfying a certain regularity condition in type $A$. Several important modules are known to be real, for example, Kirillov–Reshetikhin modules for all types \cite[Theorem 3.1]{HL16} (by a result from \cite{Qin17}), snake modules in types $A$ and $B$ \cite{DLL19,DS20}, and Hernandez--Leclerc modules \cite{BC19,DS25,GDL22,Qin17}. One of this paper's main results addresses Problem \ref{Hernandez and Leclerc's classification open question} for reachable real modules, where reachable real modules in a monoidal categorification of a cluster algebra are those modules obtained from an initial monoidal cluster via mutations.

Recall that \textit{reachable rigid objects} in an additive categorification of a cluster algebra are those objects obtained from an initial cluster-tilting object via categorical mutations. Qin in \cite{Qin20} posed \textit{the multiplicative reachability conjecture}: the real objects are reachable. It was proved in \cite[Corollary 8.6]{GLS13} that the module corresponding to a cluster variable in $\mathscr{C}^{\leq \xi}_l$ is prime. So a bijection exists between the indecomposable reachable rigid objects in the cluster category $\mathcal{C}^{\leq \xi}_l$ associated with the principal quiver of $\Gamma^{\leq \xi}_l$ and the non-frozen real prime objects in $\mathscr{C}^{\leq \xi}_l$ if the multiplicative reachability conjecture holds.

Inspired by the work of these pioneers, we formulate the following conjecture.

\begin{conjecture}[Conjecture \ref{conjecture on real prime modules}]\label{conjecture 1.2}
A non-frozen real prime module in $\mathscr{C}^{\leq \xi}_l$ is uniquely determined by the extended $\mathbf{g}$-vector of its corresponding indecomposable reachable rigid object in $\mathcal{C}^{\leq \xi}_l$ up to equivalence. Specifically, the map
\begin{align*}
\Phi^{\leq \xi}_{l}: \left\{ \makecell{\text{indecomposable reachable} \\ \text{rigid objects in $\mathcal{C}^{\leq \xi}_l$}} \right\} & \longrightarrow  \left\{ \makecell{\text{the classes of non-frozen}\\ \text{real prime objects in $\mathscr{C}^{\leq \xi}_l$}} \right\} \\
M & \longmapsto \left[L\left( (\mathbf{z}^{\leq \xi}_l)^{\widetilde{\mathbf{g}}(M)} \right)\right]
\end{align*}
defines a bijection. Here $\widetilde{\mathbf{g}}(M)$ denotes the extended $\mathbf{g}$-vector of $M$, and
\[
(\mathbf{z}^{\leq \xi}_l)^{\widetilde{\mathbf{g}}(M)} = \frac{ (\mathbf{z}^{\leq \xi}_{l-1})^{\mathbf{g}(M)} }{ F_M|_{\mathbb{P}}\bigl( (y_{i,r})_{(i,r)\in \widehat{I}^{\leq \xi}_{l-1}} \bigr) },
\]
where $F_M$ is the $F$-polynomial of $M$, and  $\mathbb{P}$ with variables $y_{i,r}$ are defined in Section \ref{Hernandez--Leclerc's classification problem}.
\end{conjecture}

For $l=1$, the map $\Phi^{\leq \xi}_1$ coincides with the bijection in \cite[Section 5.2]{DS25}, and hence the images of the indecomposable reachable rigid objects in $\mathcal{C}^{\leq \xi}_l$ under $\Phi^{\leq \xi}_{l}$ recover the Hernandez--Leclerc modules \cite{BC19,DS25}. In this case, $\widetilde{\mathbf{g}}(M)$ admits a quiver-theoretic interpretation:
\[
\widetilde{\mathbf{g}}(M) =
\begin{pmatrix}
\mathbf{g}(M) \\
\underline{\dim}(\operatorname{soc}(\mathscr{F}_{Q_{\xi}} M))
\end{pmatrix},
\]
where $\mathbf{g}(M)$ is the $\mathbf{g}$-vector of $M \in \mathcal{C}^{\leq \xi}_1$, $Q_{\xi}$ is the Dynkin quiver associated to $\xi$, and $\mathscr{F}_{Q_{\xi}}$ is the functor from $\mathcal{C}^{\leq \xi}_1$ to the category of finitely generated $\mathbb{K}Q_{\xi}$-modules defined in \cite{BMR07}.

\begin{remark}
For $l=\infty$, the quiver $\Gamma^{\leq \xi}_l$ has no frozen vertices and $\mathcal{C}^{\leq \xi}_l$ is the cluster category of $\Gamma^{\leq \xi}_l$. In this case, $\Phi^{\leq \xi}_{l}$ induces the following bijection:
\begin{align*}
\left\{ \makecell{\text{indecomposable reachable} \\ \text{rigid objects in $\mathcal{C}^{\leq \xi}_l$}} \right\} & \longrightarrow  \left\{ \makecell{\text{the classes of real prime objects in $\mathscr{C}^{\leq \xi}_l$}} \right\} \\
M & \longmapsto  \left[L\left( (\mathbf{z}^{\leq \xi}_l)^{\mathbf{g}(M)}\right)\right],
\end{align*}
where $\mathbf{g}(M)$ is the $\mathbf{g}$-vector of $M \in \mathcal{C}^{\leq \xi}_l$.
\end{remark}

We prove the following theorem.

\begin{theorem}[Theorem \ref{the classification of real modules under reachability}]\label{theorem 1.3}
Assume that the multiplicative reachability conjecture holds. Then Conjecture \ref{conjecture 1.2} holds.
\end{theorem}

We extend the map $\Phi^{\leq \xi}_{l}$ to decomposable reachable rigid objects $M_1 \oplus M_2$ with two indecomposable summands by the formula
\[
\Phi^{\leq \xi}_{l}(M_1 \oplus M_2) = \left[\Phi^{\leq \xi}_{l}(M_1)\right] \left[\Phi^{\leq \xi}_{l}(M_2) \right].
\]

By combining this with the proof of Theorem \ref{theorem 1.3}, we have the following corollary.

\begin{corollary}[Corollary \ref{induce Hernandez and Leclerc's classification}]
The map $\Phi^{\leq \xi}_{l}$ induces a bijection from the set of reachable rigid objects in $\mathcal{C}^{\leq \xi}_l$ to the set of classes of the non-frozen reachable real objects in $\mathscr{C}^{\leq \xi}_l$. Moreover, if the multiplicative reachability conjecture holds, then $\Phi^{\leq \xi}_{l}$ induces a bijection from the set of reachable rigid objects in $\mathcal{C}^{\leq \xi}_l$ to the set of classes of the non-frozen real objects in $\mathscr{C}^{\leq \xi}_l$.
\end{corollary}

Since the cluster category $\mathcal{C}^{\leq \xi}_l$ is very large and complicated in general, we consider the cluster category $\mathcal{C}_{Q_{\xi}}$  associated to a Dynkin quiver $Q_{\xi}$ instead. For an indecomposable object $M\in \mathcal{C}_{Q_{\xi}}$ and an integer $l \in \mathbb{Z}_{\geq 1}$, we obtain the following map
\begin{align*}
\Psi^{\leq \xi}_l: \left\{ \makecell{\text{rigid objects in $\mathcal{C}_{Q_{\xi}}$}} \right\} & \longrightarrow  \left\{ \makecell{ \text{simple objects in $\mathscr{C}^{\leq \xi}_l$}} \right\},
\end{align*}
see Definition \ref{a new family of real prime modules}. We call the image of a rigid object under the map $\Psi^{\leq \xi}_l$ a \textit{cluster module}. 

\begin{theorem}[Theorem \ref{indecomposable are real prime}, Corollary \ref{HL are cluster modules} and Theorem \ref{cluster modules are reachable}]
\begin{itemize}
\item[(1)] The map $\Psi^{\leq \xi}_l$ reduces to the following map 
\begin{align*}
\Psi^{\leq \xi}_l: \left\{ \makecell{\text{indecomposable objects in $\mathcal{C}_{Q_{\xi}}$}} \right\} & \longrightarrow  \left\{ \makecell{ \text{real prime modules in $\mathscr{C}^{\leq \xi}_l$}} \right\}.
\end{align*}
\item[(2)] Every Hernandez--Leclerc module is a prime cluster module.
\item[(3)] The cluster modules are reachable and real.
\end{itemize}
\end{theorem}


By employing Proposition \ref{prop 4.6} and Theorem \ref{all cluster relations associated to quiver I}, we establish the following result.

\begin{theorem}[Theorem \ref{cluster exchange relations}] \label{theorem 1.7}
Let $(L,N)$ be an exchange pair in $\mathcal{C}_{Q_{\xi}}$ with exchange triangles (\ref{non-split triangles}). Then the following equation holds in the Grothendieck ring $K_0(\mathscr{C}^{\leq \xi}_l)$:
\begin{align*}
\begin{split}
\left[ \Psi^{\leq \xi}_{l}(L) \right] \left[ \Psi^{\leq \xi}_{l}(N) \right] & = \left[ \Psi^{\leq \xi}_{l}(M) \right] \left(\prod_{i \in I} \left[ L(u_i(l)) \right]^{\alpha_i} \left[ L(v_i(l)) \right]^{\alpha_i} \right) \\
& \quad + \left[ \Psi^{\leq \xi}_{l}(M') \right] \left( \prod_{j \in I} \left[ L(u_j(l)) \right]^{\beta_j} \left[ L(v_j(l)) \right]^{\beta_j} \right),
\end{split}
\end{align*}
where the vectors $(\alpha_i)_{i\in I}$ and $(\beta_i)_{i\in I}$ are given by
\begin{align*}
(\alpha_i)_{i\in I} & = \kappa (L, M, N), \\
(\beta_j)_{j\in I} & = \kappa (L, M', N) + \mathbf{g}(\operatorname{Im}(h)),
\end{align*}
see (\ref{scalar kappa}). Moreover, the highest $\ell$-weight monomial of $\Psi^{\leq \xi}_{l}(L) \otimes \Psi^{\leq \xi}_{l}(N)$ coincides with the highest $\ell$-weight monomial of $\Psi^{\leq \xi}_{l}(M) \otimes \left( \bigotimes_{i \in I} \left( L(u_i(l)) \otimes L(v_i(l)) \right)^{\otimes \alpha_i} \right)$.
\end{theorem}

This system of equations in Theorem \ref{theorem 1.7} generalizes the relations for Hernandez--Leclerc modules established in \cite[Theorem 5.8]{DS25}, as well as the classical T-system relations for Kirillov--Reshetikhin modules, see Theorem \ref{KR are Cluster module}. Thus the cluster modules generalizes Hernandez--Leclerc modules and Kirillov--Reshetikhin modules.

This paper is organized as follows. Section \ref{preliminaries} provides necessary background on cluster algebras, quantum affine algebras, and their representations. In Section \ref{on the classification of the real prime modules}, we recall the cluster structure on $K_0(\mathscr{C}^{\leq \xi}_l)$ and state Conjecture \ref{conjecture on real prime modules}. We verify the conjecture under a multiplicative reachability condition, answering the Hernandez--Leclerc classification problem in monoidal categorifications of cluster algebras under this condition. In Section \ref{definition of cluster modules}, we introduce and study cluster modules. In particular, the Kirillov--Reshetikhin modules and Hernandez--Leclerc modules are cluster modules. In Section \ref{a new family of real prime modules and Equational systems}, we construct a cluster algebra with coefficients, in particular, we determine its extended $\mathbf{g}$-vectors in Proposition \ref{an expression of extended g-vectors}, and a system of equations in Theorem \ref{all cluster relations associated to quiver I}. These results will be used to prove Theorem \ref{cluster exchange relations}, which naturally generalizes the classical T-system relations for Kirillov--Reshetikhin modules and the exchange relations for Hernandez--Leclerc modules. In Section \ref{Proofs of two theorems}, we prove Theorem \ref{indecomposable are real prime} and Theorem \ref{KR are Cluster module}.

\textbf{Conventions.} Throughout this work, $\mathbb{K}$ denotes a field of characteristic $0$. We adopt the notation $[x]_+ = \max(x, 0)$ for integers $x \in \mathbb{Z}$. All cluster algebras are assumed to be skew-symmetric, furthermore, we identify any skew-symmetric integer matrix with its corresponding quiver.

\section{Cluster algebras, quantum affine algebras and their representations}\label{preliminaries}

\subsection{Cluster algebras}
Recall that $\mathbb{P}=\text{Trop}(u_1,\ldots,u_m)$ is a \textit{tropical semifield} if $\mathbb{P}$ is an abelian multiplicative group freely generated by the elements $u_1,\ldots,u_m$, with auxiliary addition $\oplus$ given by
\[
\prod_{j=1}^{m} u^{a_j}_j  \oplus  \prod_{j=1}^{m} u^{b_j}_j  =  \prod_{j=1}^{m} u^{\min\{a_j,b_j\}}_j.
\]
Let $\mathcal{F}=\mathbb{QP}(v_1,\ldots,v_n)$ be the field of rational functions in $n$ independent variables $v_1, \ldots, v_n$ with coefficients in $\mathbb{QP}$.

\textit{A labeled seed} $\Sigma$ in $\mathcal{F}$ is a triple $(\mathbf{x},\mathbf{y},B)$, where
\begin{itemize}
\item the \textit{cluster} $\mathbf{x}=(x_1,\ldots,x_n)$ is an $n$-tuple of cluster variables in $\mathcal{F}$ forming a free generating set over $\mathbb{QP}$,
\item the \textit{coefficient tuple} $\mathbf{y}=(y_1,\ldots,y_n)$ is an $n$-tuple of elements of $\mathbb{P}$, and
\item the \textit{exchange matrix} $B=(b_{ij})$ is an $n\times n$ skew-symmetric integer matrix, i.e., $B=-B^{T}$, where $B^T$ is the transpose of $B$.
\end{itemize}

Let $\Sigma=(\mathbf{x},\mathbf{y},B)$ be a labeled seed. For any $k\in \{1,\ldots,n\}$, a \textit{labeled seed mutation} $\mu_{k}$ in direction $k$ transforms $\Sigma$ into a new labeled seed $\mu_k(\Sigma)=\Sigma'=(\mathbf{x}',\mathbf{y}',B')$ by the following component-wise transformations:
\begin{itemize}
\item The new cluster $\mathbf{x}'=(x'_1,\ldots,x'_n)$ is given by $x'_j=x_j$ for $j\neq k$ and the following relation
\begin{align*}
x'_k = \cfrac{y_k \prod_{i=1}^{n} x^{[b_{ik}]_+}_i +  \prod_{i=1}^{n} x^{[-b_{ik}]_+}_i  }{(y_k\oplus 1)x_k}.
\end{align*}
\item The new coefficient tuple $\mathbf{y}'=(y'_1,\ldots,y'_n)$ is given by $y'_k=y^{-1}_k$ and for $j \neq k$
\begin{align*}
y'_j = y_j y^{[b_{kj}]_+}_k (y_k \oplus 1)^{-b_{kj}}.
\end{align*}
\item The new exchange matrix $B'=(b'_{ij})$ is given by
\begin{align*}
b'_{ij} = \begin{cases}
-b_{ij}  & \text{if $i=k$ or $j=k$}, \\
b_{ij} + \cfrac{|b_{ik}|b_{kj}+b_{ik}|b_{kj}|}{2}  & \text{otherwise}.
\end{cases}
\end{align*}
\end{itemize}

Let $\mathbb{T}_n$ be an $n$-regular tree whose edges are labeled by $1,\ldots,n$ so that the $n$ edges emanating from each vertex have different labels. A \textit{cluster pattern} on $\mathbb{T}_n$ is an assignment of a labeled seed $(\textbf{x}_t,\textbf{y}_t,B_t)$ to every vertex $t\in \mathbb{T}_n$ such that two labeled seeds assigned to the endpoints of an edge $t-t'$ labeled by $k$ are obtained from each other by labeled seed mutation in direction $k$. For convenience, fix a vertex $t_0\in \mathbb{T}_n$, the cluster pattern on $\mathbb{T}_n$ is uniquely determined by a labeled seed at $t_0$, which can be chosen arbitrarily. An initial labeled seed at $t_0$ is written as $\Sigma_{t_0}=(\textbf{x}_{t_0},\textbf{y}_{t_0},B_{t_0})$, where
\begin{align}\label{an initial labeled seed}
\textbf{x}_{t_0} = (x_{1;t_0},\ldots,x_{n;t_0}), \quad \textbf{y}_{t_0}=(y_{1;t_0},\ldots,y_{n;t_0}), \quad B_{t_0}=(b^{t_0}_{ij}).
\end{align}

Let $\mathcal{X}$ be the union of all cluster variables obtained from the initial labelled seed $\Sigma_{t_0}$ via mutations. The cluster algebra $\mathcal{A}(\Sigma_{t_0})$ with initial labelled seed $\Sigma_{t_0}$ is the $\mathbb{Z}\mathbb{P}$-subalgebra of the ambient field $\mathcal{F}$ generated by $\mathcal{X}$.

\subsection{Principal coefficients, $F$-polynomials, $\mathbf{g}$-vectors and $\mathbf{c}$-vectors}

A cluster pattern $t\mapsto (\textbf{x}_t,\textbf{y}_t,B_t)$, or the corresponding cluster algebra $\mathcal{A}$, has \textit{principal coefficients} at a vertex $t_0$ if $\mathbb{P}=\text{Trop}(y_1,\ldots,y_n)$ and $\textbf{y}_{t_0}=(y_1,\ldots,y_n)$. In this case, we denote by $\mathcal{A}_{\bullet}(B_{t_0})$ the cluster algebra with principal coefficients. Every cluster variable $x_{j;t}$ in $\mathcal{A}_{\bullet}(B_{t_0})$ is a Laurent polynomial $X_{j;t}(x_{1;t_0},\ldots,x_{n;t_0})$ with coefficients in $\mathbb{Z}_{>0}[y_{1;t_0},\ldots,y_{n;t_0}]$ by the Laurent property \cite[Proposition 3.6]{FZ07} and a result from \cite{LS15}.

The \textit{$F$-polynomial} $F^{B_{t_0};t_0}_{j;t}$ of $x_{j;t}$ is defined as a specialization of $X_{j;t}(x_{1;t_0},\ldots,x_{n;t_0})$ by setting all $x_{i;t_0}$ to 1. It was conjectured in \cite[Conjecture 5.4]{FZ07} and proved in \cite{DWZ10} that each $F$-polynomial $F^{B_{t_0};t_0}_{j;t}$ has constant term 1, using quivers with potentials and their representations.

For $\mathcal{A}_{\bullet}(B_{t_0})$, Fomin and Zelevinsky \cite{FZ07} introduced a $\mathbb{Z}^n$-grading structure on the Laurent polynomial ring $\mathbb{Z}[x^{\pm 1}_{1;t_0}, \ldots, x^{\pm 1}_{n;t_0}, y_1, \ldots, y_n]$ by setting
\[
\text{deg}(x_{i;t_0}) = \textbf{e}_i, \quad \text{deg}(y_j) = -\textbf{b}_j,
\]
where $i,j \in \{1,2,\ldots,n\}$, $\textbf{e}_i$ is the $i$-th column vector of the identity matrix $I_n$, and $\textbf{b}_j$ is the $j$-th column vector of $B_{t_0}$. Under this grading, they proved in \cite[Proposition 6.1]{FZ07} that every cluster variable in $\mathcal{A}_{\bullet}(B_{t_0})$ is homogeneous, and the homogeneous degree of $x_{j;t}$ is called the \textit{$\mathbf{g}$-vector} of $x_{j;t}$, denoted by $\mathbf{g}(x_{j;t})$.

The coefficient tuple $\textbf{y}_t$ and the exchange matrix $B_{t}$ can be encoded into a larger $2n\times n$ matrix $\widetilde{B}_t$ with the top $n\times n$ block $B_t$ and the bottom block $C^{B_{t_0}}_t$. The matrix $C^{B_{t_0}}_t$ is called a $C$-matrix of $\textbf{x}_t$ with the columns $c(y_{1;t}),\ldots,c(y_{n;t})$ called \textit{$\mathbf{c}$-vectors}. Moreover, $\mathbf{c}$-vectors have the following crucial property (also called the sign-coherence property of $\mathbf{c}$-vectors, see \cite{GHKK18,NZ12})
\[
\text{Each column vector of $C$-matrix is either all non-negative or all non-positive}.
\]

With the help of \cite[Proposition 1.3]{NZ12}, for an edge $t-t'$ labeled by $k$, we have the following equation.
\begin{align*}
\mathbf{g}(x_{k;t}) + \mathbf{g}(x_{k;t'}) = \begin{cases}
\sum\limits_{j=1}^{n} [b^t_{jk}]_+ \mathbf{g}(x_{j;t})  & \text{if $\varepsilon_k(C)=-1$}, \\
\sum\limits_{j=1}^{n}  [-b^t_{jk}]_+ \mathbf{g}(x_{j;t}) & \text{if $\varepsilon_k(C)=1$}.
\end{cases}
\end{align*}

If $F$ is a subtraction-free rational expression over $\mathbb{Q}$ in $l$ variables, $\mathbb{P}$ is a semifield, and $y_1, \ldots, y_l$ are some elements of $\mathbb{P}$, then we denote by $F|_{\mathbb{P}}(y_1, \ldots, y_l)$ the evaluation of $F$ at $y_1, \ldots, y_l$.

The following proposition is crucial for us.

\begin{proposition}[{\cite[Corollary 6.3]{FZ07}}] \label{Fomin-Zelevinsky separation}
Let $\mathcal{A}$ be a cluster algebra over an arbitrary semifield $\mathbb{P}$, with a seed at an initial vertex $t_0$ given by (\ref{an initial labeled seed}). Then the cluster variables in $\mathcal{A}$ can be expressed as follows:
\begin{align} \label{separate formula for a cluster variable}
x_{l;t} = \frac{\textbf{x}^{\mathbf{g}(x_{l;t})}_{t_0}}{F^{B_{t_0};t_0}_{l;t}\mid_{\mathbb{P}}(y_1,\ldots,y_n)} F^{B_{t_0};t_0}_{l;t}\mid_{\mathcal{F}}(\widehat{y}_1,\ldots,\widehat{y}_n),
\end{align}
where $\widehat{y}_j=y_j\left(\prod_{i=1}^{n} x^{b^{t_0}_{ij}}_{i;t_0}\right)$ for any $j\in \{1,\ldots,n\}$.
\end{proposition}

Let $\widetilde{\textbf{x}}_{t_0}=(x_{1;t_0},\ldots,x_{n;t_0},u_1,\ldots,u_m)$ and
\begin{align*}
\widetilde{\textbf{x}}^{\widetilde{\mathbf{g}}(x_{l;t})}_{t_0} = \cfrac{\textbf{x}^{\mathbf{g}(x_{l;t})}_{t_0}}{F^{B_{t_0};t_0}_{l;t}\mid_{\mathbb{P}}(y_1,\ldots,y_n)}.
\end{align*}
Here $\widetilde{\mathbf{g}}(x_{l;t})\in \mathbb{Z}^{n+m}$ is called the \textit{extended $\mathbf{g}$-vector} of $x_{l;t}$. Using Equation (\ref{separate formula for a cluster variable}), we have
\begin{align*} 
x_{l;t} = \widetilde{\textbf{x}}^{\widetilde{\mathbf{g}}(x_{l;t})}_{t_0} F^{B_{t_0};t_0}_{l;t}\mid_{\mathcal{F}}(\widehat{y}_1,\ldots,\widehat{y}_n).
\end{align*}

Following \cite[Theorem 1.7]{DWZ10}, the transformation of the extended $g$-vectors of cluster variables can be described as follows: for an edge $t-t'$ labeled by $k$,
\begin{align}\label{iterate formula of g-vectors}
& \widetilde{\mathbf{g}}(x_{i;t'}) =  \begin{cases}
-\widetilde{\mathbf{g}}(x_{k;t})  & \text{if $i=k$}, \\
\widetilde{\mathbf{g}}(x_{i;t}) + b^{t}_{ik}[\widetilde{\mathbf{g}}(x_{k;t})]_+ & \text{if $i\neq k$, $b_{ik} \geq 0$}, \\
\widetilde{\mathbf{g}}(x_{i;t}) + b^{t}_{ik}[-\widetilde{\mathbf{g}}(x_{k;t})]_+  & \text{ if $i\neq k$, $b_{ik} \leq 0$}.
\end{cases}
\end{align}

\subsection{Quantum affine algebras and their representations}

Let $\g$ be a simply-laced complex simple Lie algebra with Cartan matrix $C=(C_{ij})_{i,j\in I}$, where $I=\{1,2,\ldots,n\}$ is the Dynkin index set of $\g$. The \textit{quantum affine algebra} $U_q(\widehat{\g})$ is defined as the associative $\mathbb{C}$-algebra generated by $x^{\pm}_{i,r}$ ($i\in I, r\in \mathbb{Z}$), $k^{\pm 1}_i$ ($i\in I$), $h_{i,s}$ ($i\in I, s\in \mathbb{Z}\backslash\{0\}$), and $c^{\pm 1/2}$, subject to some relations, see \cite[Theorem 12.2.1]{CP94}. In this paper, we assume that the quantum parameter $q\in \mathbb{C}^\times$ is not a root of unity.

Let $\mathscr{C}_{\mathfrak{g}}$ be the category of type 1 finite-dimensional  $U_q(\widehat{\mathfrak{g}})$-modules, see \cite{CP94,FM01}. The category $\mathscr{C}_{\mathfrak{g}}$ is a tensor category due to a Hopf algebra structure of $U_q(\widehat{\mathfrak{g}})$, and it is rigid, i.e., every object $V\in \mathscr{C}_{\mathfrak{g}}$ has its left dual $\mathscr{D}(V)$ and right dual $\mathscr{D}^{-1}(V)$ in $\mathscr{C}_{\mathfrak{g}}$. Hence we have the evaluation morphisms
\[
V \otimes \mathscr{D}(V) \to \mathbb{C} \to 0,  \quad \mathscr{D}^{-1}(V) \otimes V \to \mathbb{C}  \to 0.
\]
As a monoidal category, $\mathscr{C}_{\mathfrak{g}}$ is not braided, i.e., for any two objects $V$ and $W$ in $\mathscr{C}_{\mathfrak{g}}$, $V\otimes W \not\cong W\otimes V$ in general. We denote by $K_0(\mathscr{C}_{\mathfrak{g}})$ the Grothendieck ring of $\mathscr{C}_{\mathfrak{g}}$, and by $[M]$ the equivalence class of $M\in \mathscr{C}_{\mathfrak{g}}$.

All simple modules in $\mathscr{C}_{\mathfrak{g}}$ have been classified in terms of Drinfeld polynomials \cite{CP91,CP94}. Combining this with $q$-character theory \cite{FR98}, we see that the simple modules in $\mathscr{C}_{\mathfrak{g}}$ are parameterized by the monomials in the set of infinite variables $(Y_{i,a})_{i\in I, a\in \mathbb{C}^\times}$, called \textit{dominant monomials}. For a dominant monomial $m$, we denote by $L(m)$ the corresponding simple $U_q(\widehat{\mathfrak{g}})$-module.

The theory of \textit{$q$-characters} was introduced in \cite{FR98}, defined as a map
\[
\chi_q: K_0(\mathscr{C}_{\mathfrak{g}}) \to \mathbb{Z}[Y^{\pm1}_{i,a}]_{i\in I; a\in \mathbb{C}^\times},
\]
where $\mathbb{Z}[Y^{\pm1}_{i,a}]_{i\in I; a\in \mathbb{C}^\times}$ is the ring of Laurent polynomials in infinitely many formal variables $(Y^{\pm1}_{i,a})_{i\in I; a\in \mathbb{C}^\times}$. It was proved in \cite{FM01} that the map $\chi_q$ is an injective ring morphism. Following \cite{FM01}, for $(i,a) \in  I\times \mathbb{C}^{\times}$, let
\begin{align*}
A_{i,a} = Y_{i,aq} Y_{i,aq^{-1}} \left( \prod_{j: C_{ij}=-1} Y_{j,a}^{-1} \right).
\end{align*}
For any simple module $L(m) \in \mathscr{C}_{\mathfrak{g}}$, it follows from \cite[Theorem 4.1]{FM01} that
\begin{align} \label{formal formula of q-characters}
\chi_q\left([L(m)] \right) = m \left(1+ \sum_p M_p \right),
\end{align}
where $M_p$ is a monomial in the variables $A^{-1}_{i,a}$, with $(i,a) \in  I\times \mathbb{C}^{\times}$.

\section{Classification of real modules in monoidal categorifications of cluster algebras}\label{on the classification of the real prime modules}

\subsection{Cluster structures on $K_0(\mathscr{C}^{\leq \xi}_l)$} \label{quivers and cluster algebras}

Let $\xi: I \to \mathbb{Z}$ be a height function such that $|\xi(i)-\xi(j)|=1$ if there is an edge $i\sim j$ in the Dynkin diagram of $\mathfrak{g}$. Such $\xi$ exists, since any Dynkin diagram is a tree.

Following \cite{HL16}, let $\Gamma$ be an infinite quiver with vertex set $\widehat{I} = \{(i,p) \mid i\in I, p\in \xi(i)+2\mathbb{Z}\}$ and arrows given by
\[
(i,r) \to (j,s) \text{ if and only if } C_{ij}\neq 0, s = r + C_{ij}.
\]
Note that the quiver $\Gamma$ does not depend on the choice of $\xi$ up to isomorphism.

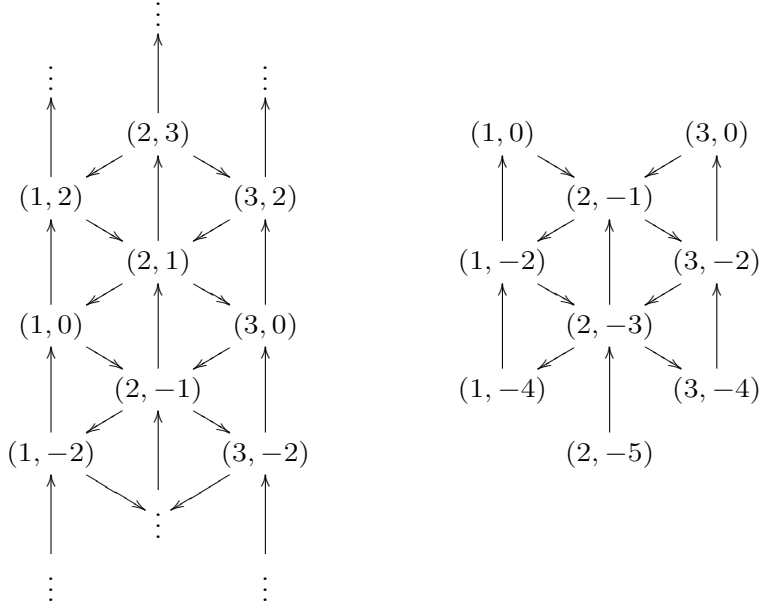
\begin{figure}[htp]
\centering
\resizebox{.65\textwidth}{.85\height}{
\begin{minipage}[c]{0.5\textwidth}
\begin{xy}
(40,60)*+{\vdots}="xx";
(25,70)*+{\vdots}="z";%
(10,60)*+{\vdots}="x";
(40,40)*+{(3,2)}="bb";
(25,50)*+{(2,3)}="a";%
(10,40)*+{(1,2)}="b";
(40,20)*+{(3,0)}="ee";
(25,30)*+{(2,1)}="d";%
(10,20)*+{(1,0)}="e";
(40,0)*+{(3,-2)}="hh";%
(25,10)*+{(2,-1)}="g";%
(10,0)*+{(1,-2)}="h";%
(10,-20)*+{\vdots}="j";
(25,-10)*+{\vdots}="k";
(40,-20)*+{\vdots}="jj";
{\ar "a";"b"};{\ar "a";"bb"};
{\ar "b";"d"};{\ar "bb";"d"};
{\ar "d";"a"};
{\ar "d";"e"};{\ar "d";"ee"};
{\ar "e";"b"};{\ar "ee";"bb"};
{\ar "e";"g"};{\ar "ee";"g"};%
{\ar "g";"d"};
{\ar "g";"h"};{\ar "g";"hh"};%
{\ar "h";"e"};{\ar "hh";"ee"};%
{\ar "j";"h"};{\ar "jj";"hh"};
{\ar "k";"g"};%
{\ar "h";"k"};{\ar "hh";"k"};%
{\ar "b";"x"};{\ar "bb";"xx"};
{\ar "a";"z"};
\end{xy}
\end{minipage}
\qquad \qquad
\begin{minipage}[c]{0.5\textwidth}
\begin{xy}
(40,40)*+{(3,0)}="bb";
(10,40)*+{(1,0)}="b";
(40,20)*+{(3,-2)}="ee";
(25,30)*+{(2,-1)}="d";%
(10,20)*+{(1,-2)}="e";
(40,0)*+{(3,-4)}="hh";%
(25,10)*+{(2,-3)}="g";%
(10,0)*+{(1,-4)}="h";%
(25,-10)*+{(2,-5)}="k";
{\ar "b";"d"};{\ar "bb";"d"};
{\ar "d";"e"};{\ar "d";"ee"};
{\ar "e";"b"};{\ar "ee";"bb"};
{\ar "e";"g"};{\ar "ee";"g"};%
{\ar "g";"d"};
{\ar "g";"h"};{\ar "g";"hh"};%
{\ar "h";"e"};{\ar "hh";"ee"};%
{\ar "k";"g"};
\end{xy}
\end{minipage}}
\caption{Quivers $\Gamma$ (left) and $\Gamma^{\leq \xi}_2$ (right) in type $A_3$, where $\xi(1)=0$, $\xi(2)=-1$ and $\xi(3)=0$.} \label{quiver Gamma in type A3}
\end{figure}

For any $l\in \mathbb{Z}_{\geq1}$, let $\Gamma^{\leq \xi}_l$ be a subquiver of  $\Gamma$ with vertex set
\[
\widehat{I}^{\leq \xi}_l = \{(i,p) \in \widehat{I} \mid \xi(i)-2l \leq p \leq \xi(i) \},
\]
where $\{ (i,\xi(i)-2l) \mid i\in I \}$ is the set of frozen vertices. It is always assumed that frozen variables cannot be mutated. Examples of quivers $\Gamma$ and $\Gamma^{\leq \xi}_2$ are shown in Figure \ref{quiver Gamma in type A3}.

For $(i,p)\in \widehat{I}^{\leq \xi}_l$, let
\begin{align}\label{the definition of zip}
z_{i,p} = \prod_{k\geq 0, p+2k \leq \xi(i)} Y_{i,p+2k}.
\end{align}
There is a cluster algebra $\mathcal{A}(\Gamma^{\leq \xi}_l)$ with initial cluster variables $\textbf{z}^{\leq \xi}_l = (z_{i,p})_{(i,p)\in \widehat{I}^{\leq \xi}_l}$.

Fix an $a\in \mathbb{C}^\times$, $i\in I$ and $r\in \mathbb{Z}$, we write $Y_{i,r}$ for $Y_{i,aq^r}$ and $A_{i,r}$ for $A_{i,aq^r}$. For $l\in \mathbb{Z}_{\geq 1}$, let
\begin{align*}
\mathcal{M} & = \left\{  \text{Monomials in the variables $(Y_{i,r})_{(i,r)\in \widehat{I}}$}   \right\}, \\
\mathcal{M}^{\leq \xi}_l & = \left\{  \text{Monomials in the variables $(Y_{i,r})_{(i,r)\in \widehat{I}^{\leq \xi}_l}$}  \right\}.
\end{align*}

\begin{definition}[{\cite[Section 3.7 and Definition 3.1]{HL10}}]  \label{the definition of subcategories}
Let $\mathscr{C}_{\mathbb{Z}}$ be the full subcategory of $\mathscr{C}_{\mathfrak{g}}$ whose objects have composition factors of the form $L(m)$, where $m \in \mathcal{M}$.

Let $\mathscr{C}^{\leq \xi}_l$ be the subcategory of $\mathscr{C}_{\mathbb{Z}}$ whose objects have composition factors of the form $L(m)$, where $m\in \mathcal{M}^{\leq \xi}_l$.
\end{definition}

Following \cite{HL10}, there is a series of monoidal subcategories of $\mathscr{C}_{\mathbb{Z}}$:
\begin{align*}
\mathscr{C}^{\leq \xi}_1 \subset \mathscr{C}^{\leq \xi}_2 \subset \cdots \subset \mathscr{C}^{\leq \xi}_l \subset \cdots \subset \cdots.
\end{align*}
Using our notation, the $q$-character map is an injective ring morphism:
\begin{align*}
\chi_q: K_0(\mathscr{C}_{\mathbb{Z}}) \to \mathbb{Z}[Y^{\pm 1}_{i,r} \mid (i,r) \in \hat{I}].
\end{align*}

Following \cite{HL10}, for $m \in \mathcal{M}^{\leq \xi}_l$, let $\chi^{\leq \xi}_q([L(m)])$ be the Laurent polynomial obtained from $\chi_q([L(m)])$ by removing some monomials in variables $Y_{i,r}$, with $(i,r)\not \in \widehat{I}^{\leq \xi}_l$. So there is an injective ring morphism:
\begin{align*}
\chi^{\leq \xi}_q:  K_0(\mathscr{C}^{\leq \xi}_l) \to \mathbb{Z}[Y^{\pm 1}_{i,r} \mid (i,r) \in \widehat{I}^{\leq \xi}_l].
\end{align*}

For $i \in I$, $k\geq 1$, $r\in \mathbb{Z}$, let
\[
m = \prod_{j=0}^{k-1} Y_{i,r+2j}.
\]
Then $W^{(i)}_{k,r}=L(m)$ is called a \textit{Kirillov--Reshetikhin module}. It was proved in \cite{Nak03}        that the $q$-characters of the Kirillov--Reshetikhin modules $W^{(i)}_{k,r}$ ($i \in I$, $k\geq 1$, $r\in \mathbb{Z}$) satisfy the corresponding T-system relation, that is,
\begin{align} \label{equ:T-system equation}
\left[ W^{(i)}_{k,r+1} \right]  \left[ W^{(i)}_{k,r-1} \right]  =  \left[ W^{(i)}_{k-1,r+1} \right]  \left[ W^{(i)}_{k+1,r-1} \right]  + \prod_{j: C_{ij}=-1} \left[ W^{(j)}_{k,r} \right],
\end{align}
where $i\in I$, $k\geq 1$, $r\in \mathbb{Z}$ and $T^{(i)}_{0,r}=1$.

\begin{theorem}[{\cite[Theorem 5.1]{HL16},\cite[Theorem 5.1]{HL21},\cite[Theorem 5.2.4]{Bit21},\cite[Theorem 4.3]{DS25}}]\label{cluster structures on C}
There is a cluster algebra isomorphism
\[
\chi_q^{\leq \xi}: K_0(\mathscr{C}^{\leq \xi}_l) \overset{\sim}{\longrightarrow} \mathcal{A}(\Gamma^{\leq \xi}_l)
\]
such that
\begin{align*}
\left[ W^{(i)}_{k+1,\xi(i)-2k}\right] \mapsto  z_{i,\xi(i)-2k}.
\end{align*}
\end{theorem}

As a cluster algebra, $ K_0(\mathscr{C}^{\leq \xi}_l)$ has frozen vertices $\left(\left[W^{(i)}_{l+1,\xi(i)-2l}\right]\right)_{i\in I}$.

\subsection{Hernandez--Leclerc's classification problem}\label{Hernandez--Leclerc's classification problem}

Let $\mathcal{Q}^{\leq \xi}_l$ be the principal quiver of $\Gamma^{\leq \xi}_l$ with vertex set $\widehat{I}^{\leq \xi}_{l-1}$. We denote by $\mathcal{C}^{\leq \xi}_l$ the associated cluster category of $\mathcal{Q}^{\leq \xi}_l$.

Let
\begin{align*}
y_{i,r} = \begin{cases}
1 & \text{if $(i,r)\in \widehat{I}^{\leq \xi}_{l-2}$}, \\
z^{-1}_{i,r-2} \prod_{(j,s)\to (i,r)} z_{j,s-2} & \text{if $(i,r)\in \widehat{I}^{\leq \xi}_{l-1} \setminus \widehat{I}^{\leq \xi}_{l-2}$}, \\
\end{cases}
\end{align*}
and $\mathbb{P}$ is the tropical semifield generated by $z_{i,r}$, with $(i,r)\in \widehat{I}^{\leq \xi}_l \setminus \widehat{I}^{\leq \xi}_{l-1}$. It is clear that $y_{i,r} \in \mathbb{P}$.  For a vertex $(i,r)$ of $\mathcal{Q}^{\leq \xi}_l$, we define $\widehat{y}_{i,r}=A^{-1}_{i,r-1}$.

\textit{Reachable real objects} (respectively, \textit{reachable real prime objects}) in a monoidal categorification of a cluster algebra are those objects obtained from an initial monoidal cluster (respectively, an arbitrary real prime object in an initial monoidal cluster) via mutations. It was conjectured in \cite{HL10,HL16,HL21,DS25,Qin20} that the real objects in $\mathscr{C}^{\leq \xi}_l$ are reachable. The correspondence above is called \textit{the multiplicative reachability conjecture} by Qin in \cite{Qin20}.

\textit{Reachable rigid objects} (respectively, \textit{indecomposable reachable rigid objects}) in an additive categorification of a cluster algebra are those objects obtained from an initial cluster-tilting object  (respectively, an initial indecomposable cluster-tilting object) via categorical mutations. Qin in \cite{Qin20} posed \textit{the additive reachability conjecture}: the rigid objects are reachable. It was proved in \cite[Corollary 8.6]{GLS13} that the cluster variables in $\mathcal{A}(\Gamma^{\leq \xi}_l)$ are the classes of real prime objects in $\mathscr{C}^{\leq \xi}_l$.

Inspired by this pioneering work, we formulate the following conjecture.

\begin{conjecture}\label{conjecture on real prime modules}
A non-frozen real prime module in $\mathscr{C}^{\leq \xi}_l$ is completely determined by the extended $\mathbf{g}$-vector of the corresponding indecomposable reachable rigid object in $\mathcal{C}^{\leq \xi}_l$ up to equivalence. More precisely, we have the following bijection:
\begin{align*}
\Phi^{\leq \xi}_{l}: \left\{ \makecell{\text{indecomposable reachable} \\ \text{rigid objects in $\mathcal{C}^{\leq \xi}_l$}} \right\} & \longrightarrow  \left\{ \makecell{\text{the classes of non-frozen}\\ \text{real prime objects in $\mathscr{C}^{\leq \xi}_l$}} \right\} \\
M & \longmapsto \left[L\left((\textbf{z}^{\leq \xi}_l)^{\widetilde{\mathbf{g}}(M)}\right)\right],
\end{align*}
where $\widetilde{\mathbf{g}}(M)$ denotes the extended $\mathbf{g}$-vector of $M$ and
\[
(\textbf{z}^{\leq \xi}_l)^{\widetilde{\mathbf{g}}(M)} = \frac{(\textbf{z}^{\leq \xi}_{l-1})^{\mathbf{g}(M)}}{F_M|_{\mathbb{P}}\left((y_{i,r})_{(i,r)\in \widehat{I}^{\leq \xi}_{l-1}}\right)}.
\]
\end{conjecture}

\begin{remark}
``Non-frozen" in Conjecture \ref{conjecture on real prime modules} comes from the fact that the cluster category $\mathcal{C}^{\leq \xi}_l$ has no information about coefficients. In order to recover the frozen real prime modules as well one would need to replace $\mathcal{C}^{\leq \xi}_l$ by other category including information of coefficients. Here the Higgs category of the Ginzburg dg-algebra associated to $\Gamma^{\leq \xi}_l$ maybe a candidate, see \cite{W23}.
\end{remark}

We are now ready for the main result of this section.

\begin{theorem}\label{the classification of real modules under reachability}
Assume that the multiplicative reachability conjecture holds. Then Conjecture \ref{conjecture on real prime modules} holds.
\end{theorem}

\begin{proof}
Assume that the multiplicative reachability conjecture holds, that is, the real objects in $\mathscr{C}^{\leq \xi}_l$ are reachable. Let $\Phi^{\leq \xi}_l$ be the composition of the bijections shown in Figure \ref{the formation of our map}.

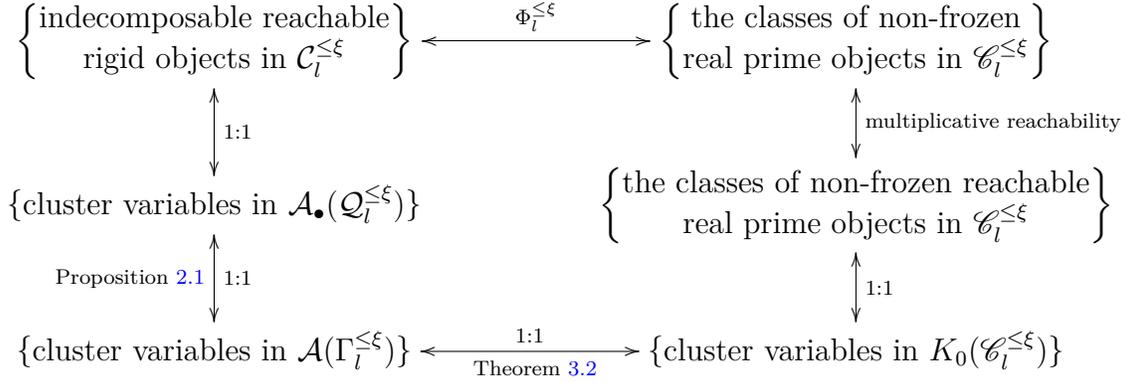
\begin{figure}
\begin{gather*}
\xymatrix@C=5em@R=5ex{ \left\{ \makecell{\text{indecomposable  reachable} \\ \text{rigid objects in $\mathcal{C}^{\leq \xi}_l$}} \right\}
\ar@{<->}[r]^-{\Phi^{\leq \xi}_l} \ar@{<->}[d]^-{1:1} & \left\{ \makecell{\text{the classes of non-frozen}\\ \text{real prime objects in $\mathscr{C}^{\leq \xi}_l$}} \right\} \\
\{ \makecell{\text{cluster variables in $\mathcal{A}_{\bullet}(\mathcal{Q}^{\leq \xi}_l)$}} \}
\ar@{<->}[d]_{\text{Proposition \ref{Fomin-Zelevinsky separation}}}^-{1:1} & \left\{ \makecell{\text{the classes of non-frozen reachable}\\ \text{real prime objects in $\mathscr{C}^{\leq \xi}_l$}} \right\}  \ar@{<->}[u]_-{\text{multiplicative reachability}} \\
\{ \makecell{\text{cluster variables in $\mathcal{A}(\Gamma^{\leq \xi}_l)$}} \} \ar@{<->}[r]_{\text{Theorem \ref{cluster structures on C}}}^-{1:1}  &
\{ \makecell{\text{cluster variables in $K_0(\mathscr{C}^{\leq \xi}_{l})$}}\} 
\ar@{<->}[u]_-{1:1}}
\end{gather*}
\caption{The formation of $\Phi^{\leq \xi}_l$.}\label{the formation of our map}
\end{figure}

By Theorem \ref{cluster structures on C}, we assume without loss of generality that $[L(m)]\in K_0(\mathscr{C}^{\leq \xi}_l)$ corresponds to a cluster variable associated to an indecomposable reachable rigid object $M \in \mathcal{C}^{\leq \xi}_l$. By Equations (\ref{separate formula for a cluster variable}) and (\ref{formal formula of q-characters}), we have
\begin{align*}
\frac{(\textbf{z}^{\leq \xi}_{l-1})^{\mathbf{g}(M)}}{F_M|_{\mathbb{P}}\left((y_{i,r})_{(i,r)\in \widehat{I}^{\leq \xi}_{l-1}}\right)} F_M\left((\widehat{y}_{i,r})_{(i,r)\in \widehat{I}^{\leq \xi}_{l-1}}\right) = m \left(1+ \sum_p M_p \right).
\end{align*}

Since $\frac{(\textbf{z}^{\leq \xi}_{l-1})^{\mathbf{g}(M)}}{F_M|_{\mathbb{P}}\left((y_{i,r})_{(i,r)\in \widehat{I}^{\leq \xi}_{l-1}}\right)}$ does not contain any $z_{i,r}$ with negative exponent, it does not contain any $Y_{i,s}$ with negative exponent. Moreover, the $F$-polynomial is a polynomial with a constant term 1. Hence
\begin{align*}
\frac{(\textbf{z}^{\leq \xi}_{l-1})^{\mathbf{g}(M)}}{F_M|_{\mathbb{P}}((y_{i,r})_{(i,r)\in \widehat{I}^{\leq \xi}_{l-1}})} = m, \quad  F_M((\widehat{y}_{i,r})_{(i,r)\in \widehat{I}^{\leq \xi}_{l-1}}) = 1+ \sum_p M_p.
\end{align*}
\end{proof}

In particular, we have the following corollary.

\begin{corollary}\label{conjecture hold for l=1}
Conjecture \ref{conjecture on real prime modules} holds for $l=1$ of all types, $l\leq 4$ of type $A_2$, and $l=2$ of type $A_3$ and type $A_4$.
\end{corollary}

\begin{proof}
For $l=1$, $l\leq 4$ of type $A_2$, and $l=2$ of type $A_3$ and type $A_4$, the involved cluster algebras are of finite type. In these cases, the multiplicative reachability conjecture holds. Therefore the corollary follows from Theorem \ref{the classification of real modules under reachability}.
\end{proof}

The real prime modules appearing in Corollary \ref{conjecture hold for l=1} have been studied and completely classified in \cite[Section 5 and Section 6.2]{DS25}. For $l=1$, these real prime modules are referred to as \textit{Hernandez--Leclerc modules}, see \cite[Section 5.2]{DS25}. We will introduce a new family of real modules, which contains Hernandez--Leclerc modules as a special subfamily in the following sections.

We extend the map $\Phi^{\leq \xi}_{l}$ to decomposable reachable rigid objects $M_1 \oplus M_2$ with two indecomposable reachable summands by the formula
\[
\Phi^{\leq \xi}_{l}(M_1 \oplus M_2) = \Phi^{\leq \xi}_{l}(M_1) \otimes \Phi^{\leq \xi}_{l}(M_2).
\]
Thus, by combining this with the proof of Theorem \ref{the classification of real modules under reachability}, we have the following corollary.

\begin{corollary}\label{induce Hernandez and Leclerc's classification}
The map $\Phi^{\leq \xi}_{l}$ induces a bijection from the set of reachable rigid objects in $\mathcal{C}^{\leq \xi}_l$ to the set of classes of the non-frozen reachable real objects in $\mathscr{C}^{\leq \xi}_l$. Moreover, if the multiplicative reachability conjecture holds, then $\Phi^{\leq \xi}_{l}$ induces a bijection from the set of reachable rigid objects in $\mathcal{C}^{\leq \xi}_l$ to the set of classes of the non-frozen real objects in $\mathscr{C}^{\leq \xi}_l$.
\end{corollary}

Since the cluster algebra $\mathcal{A}(\Gamma^{\leq \xi}_l)$ admits an additive categorification via the cluster category $\mathcal{C}^{\leq \xi}_l$, the formula (\ref{iterate formula of g-vectors}) describes the transformation rule for the extended $g$-vectors of reachable rigid objects in $\mathcal{C}^{\leq \xi}_l$. Combined with Corollary \ref{induce Hernandez and Leclerc's classification}, this allows us to address Problem \ref{Hernandez and Leclerc's classification open question} for reachable real modules. Moreover, if the multiplicative reachability conjecture holds, then Problem \ref{Hernandez and Leclerc's classification open question} can be fully solved for all real modules.

\section{Cluster modules}\label{definition of cluster modules}

Let $S \in \mathscr{C}_{\g}$ be a simple module. Denote by $\hw (S)$ the highest $\ell$-weight monomial in $\chi_q([S])$. It was shown in \cite{FM01} that, for any two simple modules $V_1$ and  $V_2$, we have
\[
\hw (V_1 \otimes V_2)=\hw (V_1) \hw(V_2).
\]
In this section, we always assume that $\xi$ is a height function, the quiver $Q_\xi$ is the Dynkin quiver associated to $\xi$, and $l\in \mathbb{Z}_{\geq 1}$.

\subsection{Cluster modules}
Recall that $z_{i,p}$, with $(i,p)\in \widehat{I}^{\leq \xi}_l$, is defined in (\ref{the definition of zip}). Let 
\[
\textbf{x}^{\leq \xi}_l =(z_{i,\xi(i)-2l+2})_{i\in I}, \quad {\bf c}^{\leq \xi}_l = (u_i(l) v_i(l))_{i\in I},   \quad  \widetilde{\textbf{x}}^{\leq \xi}_l = (\textbf{x}^{\leq \xi}_l, {\bf c}^{\leq \xi}_l),
\] 
where
\begin{align*}
u_i(l) = \cfrac{z_{i,\xi(i)-2l+2}}{z_{i,\xi(i)}},  
\quad  v_i(l) = z_{i,\xi(i)-2l}. 
\end{align*}
We denote by $\mathbb{P}_l$ the tropical semifield generated by $u_i(l)$ and $v_i(l)$, with $i\in I$.

Recall that every indecomposable object in the cluster category $\mathcal{C}_{Q_{\xi}}$ is rigid, and every rigid object is reachable.

\begin{definition}\label{a new family of real prime modules}
Let $\Psi^{\leq \xi}_l$ be a map from the set of rigid objects in $\mathcal{C}_{Q_{\xi}}$ to the set of simple objects in $\mathscr{C}^{\leq \xi}_l$, such that
\begin{itemize}
\item For an arbitrary indecomposable object $M$ in $\mathcal{C}_{Q_{\xi}}$, 
\[
\Psi^{\leq \xi}_l(M)=L\left((\widetilde{\textbf{x}}^{\leq \xi}_l)^{\widetilde{\mathbf{g}}(M)}\right),
\]
where
\[
\widetilde{\mathbf{g}}(M) =
\begin{pmatrix}
\mathbf{g}(M) \\
\underline{\dim}(\operatorname{soc}(\mathscr{F}_{Q_{\xi}} M))
\end{pmatrix}
\]
is the extended $\mathbf{g}$-vector of $M$. 
\item This map $\Psi^{\leq \xi}_l$ extends to arbitrary rigid objects in $\mathcal{C}_{Q_{\xi}}$ via the following formula:
\[
\Psi^{\leq \xi}_{l}(N_1 \oplus N_2) = \Psi^{\leq \xi}_{l}(N_1) \otimes \Psi^{\leq \xi}_{l}(N_2).
\]
\end{itemize}
\end{definition}

\begin{definition}
A simple module in $\mathscr{C}^{\leq \xi}_l$ is called a \textit{cluster module} if it is the image of a rigid object in the cluster category $\mathcal{C}_{Q_{\xi}}$ under the map $\Psi^{\leq \xi}_l$.
\end{definition}

In the following, we give an example of $\Psi^{\leq \xi}_{l}$.

\begin{example}\label{computation of the extended g-vectors}
Assume that $\xi$ is the height function in type $A_3$ such that
\[
\xi(1)=0, \,\, \xi(2)=-1 \text{ and } \xi(3)=-2.
\]
By definition, $Q_{\xi}$ is the quiver $1\to 2 \to 3$. The Auslander--Reiten quiver of the cluster category $\mathcal{C}_{Q_{\xi}}$ is shown in Figure \ref{the Auslander-Reiten quiver of linear A3}, and the corresponding extended $\mathbf{g}$-vectors are shown in Figure \ref{The extended g-vectors of linear A3}. The top part of the extended $\mathbf{g}$-vector is the usual $\mathbf{g}$-vector, and its bottom part is the dimension vector of the socle of the corresponding module. For example, the socle of the module $\substack{\,1\, \\ 2}$ is the simple module $S(2)$, and therefore, the bottom part of the extended $\mathbf{g}$-vector $\widetilde{\mathbf{g}}(\substack{\,1\, \\ 2})$ is $\left(\begin{array}{c} 0 \\ 1 \\ 0 \end{array} \right)$.

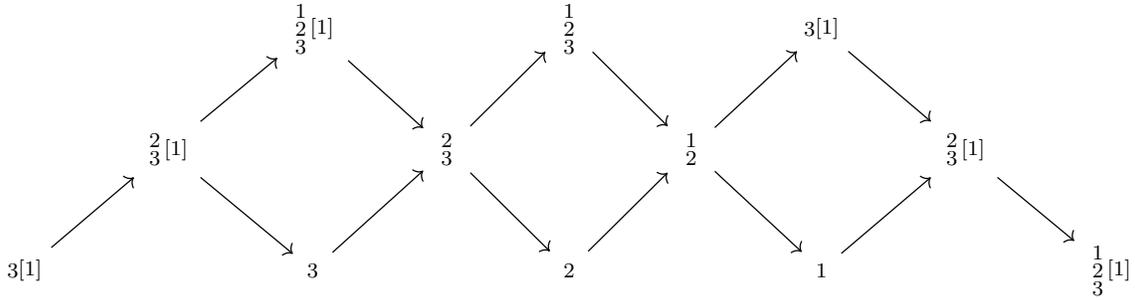
\begin{figure}[htp]
\center
\begin{tikzcd}
&& \substack{\,1 \, \\ \,2\, \\ 3} \substack{[1]} \ar[dr]  && \substack{\,1 \, \\ \,2\, \\ 3} \ar[dr] &&\substack{3}\substack{[1]} \ar[dr]  &&  \\
& \substack{\,2\, \\ 3}\substack{[1]} \ar[ur] \ar[dr]  && \substack{\,2\, \\ 3} \ar[ur] \ar[dr] &  &  \substack{\,1\, \\ 2} \ar[ur]  \ar[dr]  & &   \substack{\,2\, \\ 3}\substack{[1]} \ar[dr] & \\
\substack{3}\substack{[1]}  \ar[ur] && \substack{3} \ar[ur] && \substack{2} \ar[ur]  &&  \substack{1} \ar[ur] && \substack{\,1 \, \\ \,2\, \\ 3} \substack{[1]}
\end{tikzcd}
\caption{The Auslander--Reiten quiver of $\mathcal{C}_{Q_{\xi}}$.}\label{the Auslander-Reiten quiver of linear A3}
\end{figure}

\begin{figure}[htp]
{\tiny
\begin{tikzcd}[column sep=tiny]
&& \left( \begin{array}{c} 1 \\ 0 \\ 0 \\ \hline  0 \\ 0 \\ 0 \end{array} \right)  \ar[dr]  
&& \left( \begin{array}{c} 0 \\ 0 \\ -1 \\  \hline  0 \\ 0 \\ 1 \end{array} \right) \ar[dr]       
&& \left(\begin{array}{c} 0 \\ 0 \\ 1 \\ \hline 0 \\ 0 \\ 0 \end{array} \right)  \ar[dr]  &&  \\
& \left(\begin{array}{c} 0 \\ 1 \\ 0 \\  \hline  0 \\ 0 \\ 0 \end{array} \right)  \ar[ur] \ar[dr]  && \left(\begin{array}{c} 1 \\ 0 \\ -1 \\ \hline 0 \\ 0 \\ 1 \end{array} \right) \ar[ur] \ar[dr]  &&  \left(\begin{array}{c} 0 \\ -1 \\ 0 \\  \hline 0 \\ 1 \\ 0 \end{array} \right) \ar[ur]  \ar[dr]  & &   \left(\begin{array}{c} 0 \\ 1 \\ 0 \\ \hline 0 \\ 0 \\ 0 \end{array} \right) \ar[dr] & \\
\left(\begin{array}{c} 0 \\ 0 \\ 1 \\ \hline 0 \\ 0 \\ 0 \end{array} \right) \ar[ur] && \left(\begin{array}{c} 0 \\ 1 \\ -1 \\ \hline 0 \\ 0 \\ 1 \end{array} \right) \ar[ur] && \left(\begin{array}{c} 1 \\ -1 \\ 0 \\ \hline 0 \\ 1 \\ 0 \end{array} \right) \ar[ur]  &&  \left(\begin{array}{c} -1 \\ 0 \\ 0 \\ \hline 1 \\ 0 \\ 0 \end{array} \right) \ar[ur]  &&  
\left( \begin{array}{c} 1 \\ 0 \\ 0 \\ \hline 0 \\ 0 \\ 0 \end{array} \right)
\end{tikzcd}}
\caption{Extended $\mathbf{g}$-vectors of indecomposable objects in $\mathcal{C}_{Q_{\xi}}$.}\label{The extended g-vectors of linear A3}
\end{figure}
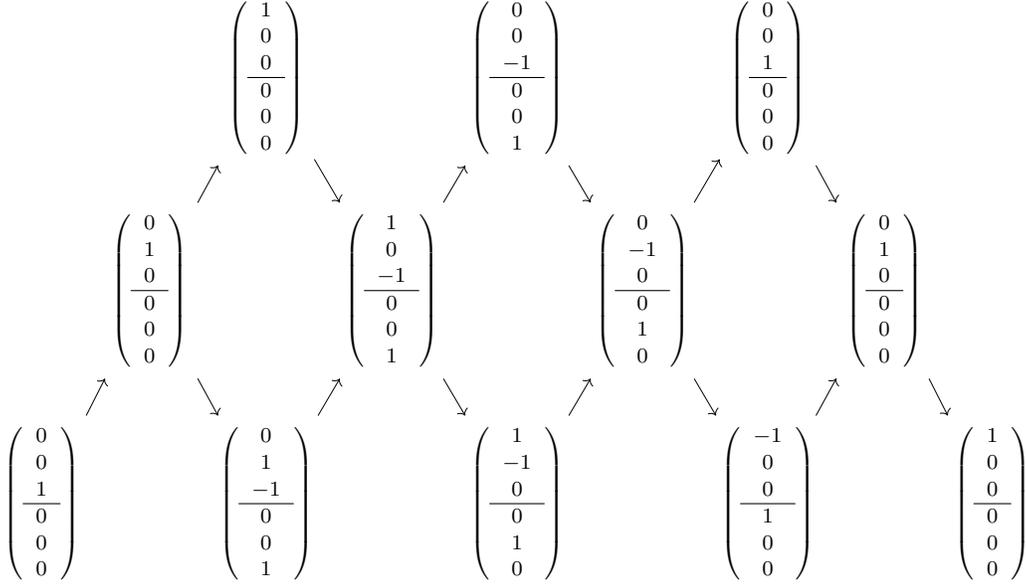


Let $l=2$. By Definition \ref{a new family of real prime modules}, the highest $\ell$-weight monomial of $\Psi^{\leq \xi}_2(\substack{\,2\, \\ 3})$ is
\[
(\widetilde{\textbf{x}}^{\leq \xi}_1)^{\widetilde{\mathbf{g}}(\substack{\,2\, \\ 3})} = z_{1,\xi(1)-2} z^{-1}_{3,\xi(3)-2} \cfrac{z_{3,\xi(3)-2}}{z_{3,\xi(i)}} z_{3,\xi(3)-4} = Y_{1,-2}Y_{1,0}  Y_{3,-6}Y_{3,-4}.
\]
These highest $\ell$-weight monomials of the images of arbitrary indecomposable objects in $\mathcal{C}_{Q_{\xi}}$ under the map $\Psi^{\leq \xi}_2$ are computed as shown in Figure \ref{the highest l-weight monomials of the prime cluster modules in level 2 of linear A3}.


\begin{figure}[H]
\[
\scalebox{0.65}{
\xymatrix@C=0.3em{
&& Y_{1,-2}Y_{1,0} \ar[dr] && Y_{3,-6}Y_{3,-4} \ar[dr] && Y_{3,-4}Y_{3,-2} \ar[dr] && \\
  & Y_{2,-3}Y_{2,-1} \ar[ur] \ar[dr] && Y_{1,-2}Y_{1,0}Y_{3,-6}Y_{3,-4} \ar[ur] \ar[dr] && Y_{2,-5}Y_{2,-3} \ar[ur] \ar[dr] && Y_{2,-3}Y_{2,-1} \ar[dr] & \\
  Y_{3,-4}Y_{3,-2} \ar[ur] && Y_{2,-3}Y_{2,-1}Y_{3,-6}Y_{3,-4} \ar[ur] && Y_{1,-2}Y_{1,0}Y_{2,-5}Y_{2,-3} \ar[ur] && Y_{1,-4}Y_{1,-2} \ar[ur] && Y_{1,-2}Y_{1,0}}
}
\]
\caption{Highest $\ell$-weight monomials of the images of arbitrary indecomposable objects in $\mathcal{C}_{Q_{\xi}}$ under the map $\Psi^{\leq \xi}_2$.}\label{the highest l-weight monomials of the prime cluster modules in level 2 of linear A3}
\end{figure}
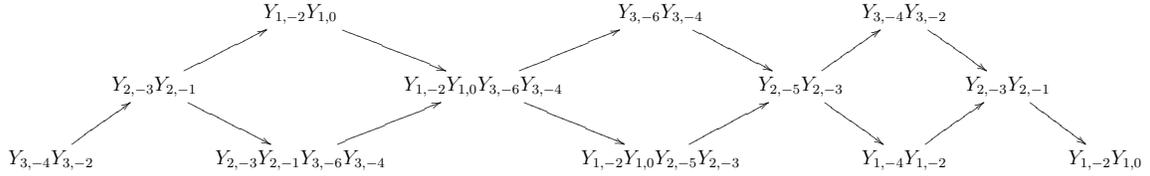
\end{example}

\begin{theorem}\label{property of our map}
\begin{itemize}
\item[(1)] The map $\Psi^{\leq \xi}_{l}$ is injective.

\item[(2)] For any indecomposable projective object $P(i)$ at $i$, with $i\in I$,
\[
\Psi^{\leq \xi}_{l}(P(i)[1]) = W^{(i)}_{l,\xi(i)-2l+2}.
\]

\item[(3)] For any indecomposable injective object $I(i)$ at $i$, with $i\in I$,
\[
\Psi^{\leq \xi}_{l}(I(i)) = W^{(i)}_{l,\xi(i)-2l}.
\]
\end{itemize}
\end{theorem}
\begin{proof}
(1) Since two non-isomorphic rigid objects in $\mathcal{C}_{Q_{\xi}}$ have different $\mathbf{g}$-vectors, by definition, the map $\Psi^{\leq \xi}_{l}$ is injective.

(2) For any $i\in I$, $\Psi^{\leq \xi}_{l}(P(i)[1])$ follows from
\[
\widetilde{\mathbf{g}}(P(i)[1]) =
\begin{pmatrix}
\mathbf{e}_i \\
\mathbf{0}
\end{pmatrix}.
\]

(3) For any $i\in I$, we have
\[
\widetilde{\mathbf{g}}(I(i)) =
\begin{pmatrix}
-\mathbf{e}_i \\
\mathbf{e}_i
\end{pmatrix}.
\]
By definition, $\Psi^{\leq \xi}_{l}(I(i)) = L\left(\cfrac{u_i(l) v_i(l)}{z_{i,\xi(i)-2l+2}}\right) = W^{(i)}_{l,\xi(i)-2l}$.
\end{proof}

\subsection{Reachability and reality}\label{A cluster subalgebra of C}
In this section, we build a one-parameter family of cluster algebras inside $K_0(\mathscr{C}^{\leq \xi}_l)$. Recall that $I=\{1,2,\ldots,n\}$ is the Dynkin index set of a simply-laced complex simple Lie algebra $\g$.

\begin{definition}\label{defintion of the quiver}
Let $Q^{\leq \xi}_l$ be a quiver with vertex set
\[
\{ (i,\xi(i)-2l+4),  (i,\xi(i)-2l+2), (i,\xi(i)-2l) \mid i\in I\},
\]
where we set $(i, r_i) = \emptyset$ when $r_i > \xi(i)$, and with arrows determined by:
\begin{enumerate}
\item[(1)] All vertices in $\{ (i,\xi(i)-2l+4), (i,\xi(i)-2l) \mid i\in I\}$ are frozen, with no arrows between any frozen vertices.
\item[(2)] For $2 \leq j \leq n$ with $\xi(j) = \xi(j-1) + 1$, the arrows incident to vertex $(j,\xi(j)-2l+2)$ are given in Figure~\ref{arrows of our quiver 1}.
\begin{figure}[htp]
\begin{align*}
& \resizebox{.8\textwidth}{.8\height}{ \xymatrix{
&  (j,\xi(j)-2l+4) \ar[d]  &   \\
(j-1,\xi(j-1)-2l+2) & (j,\xi(j)-2l+2) \ar[l] \ar[r] & (j+1,\xi(j+1)-2l+2)  \\
 & (j,\xi(j)-2l) \ar[u]  & }} \\ 
& \resizebox{.8\textwidth}{.8\height}{\xymatrix{
& (j,\xi(j)-2l+4) \ar[d]  & (j+1,\xi(j+1)-2l+4) \\
(j-1,\xi(j-1)-2l+2) & (j,\xi(j)-2l+2) \ar[l]  \ar[ur] \ar[dr]  & (j+1,\xi(j+1)-2l+2) \ar[l] \\
& (j,\xi(j)-2l) \ar[u]  & (j+1,\xi(j+1)-2l)}}
\end{align*}
\caption{Local quiver at $(j,\xi(j)-2l+2)$ when $\xi(j) = \xi(j+1)+1$ (top) and $\xi(j) = \xi(j+1)-1$ (bottom).}\label{arrows of our quiver 1}
\end{figure}
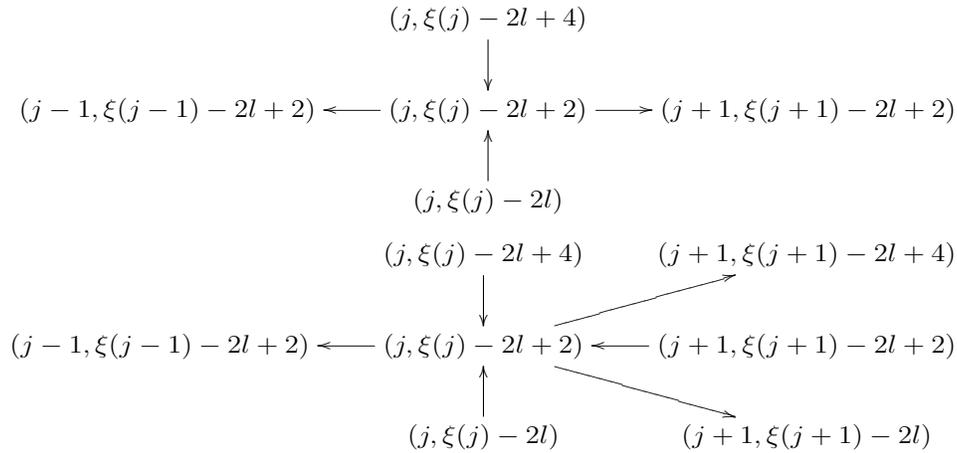
\item[(3)] For $2 \leq j \leq n$ with $\xi(j) = \xi(j-1) - 1$, the arrows incident to vertex $(j,\xi(j)-2l+2)$ are given in Figure~\ref{arrows of our quiver 2}.
\begin{figure}[htp]
\begin{align*}
& \resizebox{.8\textwidth}{.8\height}{\xymatrix{
(j-1,\xi(j-1)-2l+4) & (j,\xi(j)-2l+4) \ar[d]  &   \\
(j-1,\xi(j-1)-2l+2) \ar[r] & (j,\xi(j)-2l+2) \ar[r] \ar[ul] \ar[dl] & (j+1,\xi(j+1)-2l+2)  \\
(j-1,\xi(j-1)-2l)  & (j,\xi(j)-2l)  \ar[u]  & }} \\
& \resizebox{.8\textwidth}{.8\height}{ \xymatrix{
(j-1,\xi(j-1)-2l+4) & (j,\xi(j)-2l+4)\ar[d]  & (j+1,\xi(j+1)-2l+4) \\
(j-1,\xi(j-1)-2l+2) \ar[r] & (j,\xi(j)-2l+2)  \ar[ur] \ar[dr]  \ar[ul] \ar[dl] & (j+1,\xi(j+1)-2l+2) \ar[l] \\
(j-1,\xi(j-1)-2l) & (j,\xi(j)-2l) \ar[u]  & (j+1,\xi(j+1)-2l)}}
\end{align*}
\caption{Local quiver at $(j,\xi(j)-2l+2)$ when $\xi(j) = \xi(j+1)+1$ (top) and $\xi(j) = \xi(j+1)-1$ (bottom).} \label{arrows of our quiver 2}
\end{figure}
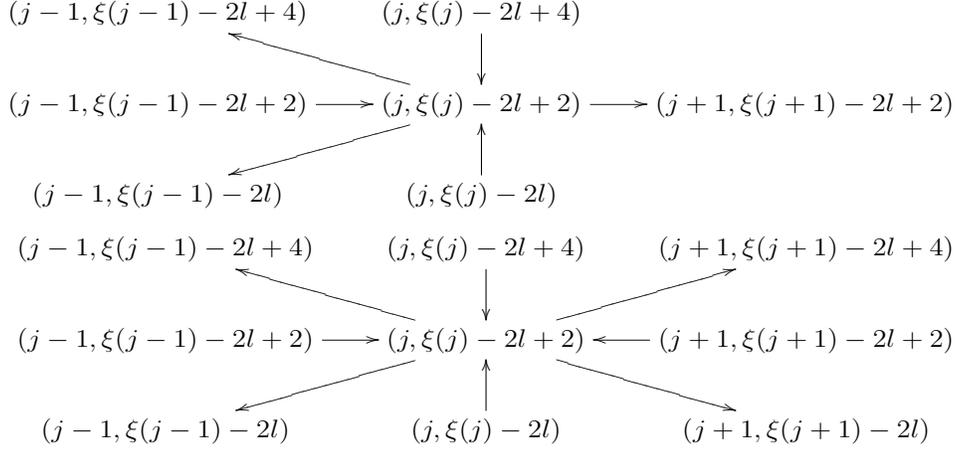
\end{enumerate}
\end{definition}

\begin{example}
Let $\g$ be of type $A_4$, and $\xi(1)=0$, $\xi(2)=-1$, $\xi(3)=-2$, $\xi(4)=-1$. By Definition \ref{defintion of the quiver}, the quiver $Q^{\leq \xi}_2$ is shown in Figure \ref{quiver in type A4}.

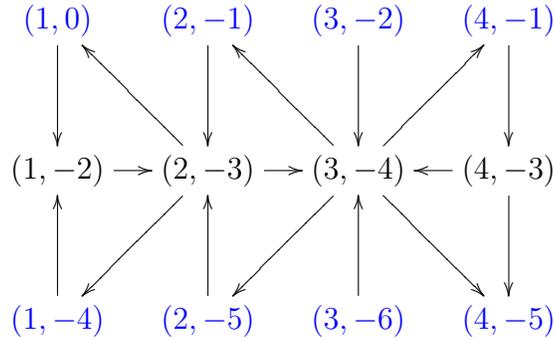
\begin{figure}[H]
\centerline{
\begin{xy}
(10,40)*+{\color{blue}{(1,0)}}="b";
(30,40)*+{\color{blue}{(2,-1)}}="d";%
(50,40)*+{\color{blue}{(3,-2)}}="ee";
(70,40)*+{\color{blue}{(4,-1)}}="ff";
(10,20)*+{(1,-2)}="e";
(30,20)*+{(2,-3)}="g";%
(50,20)*+{(3,-4)}="hh";%
(70,20)*+{(4,-3)}="hhh";%
(10,0)*+{\color{blue}{(1,-4)}}="h";%
(30,0)*+{\color{blue}{(2,-5)}}="k";
(50,0)*+{\color{blue}{(3,-6)}}="kk";%
(70,0)*+{\color{blue}{(4,-5)}}="jj";%
{\ar "g";"b"};{\ar "b";"e"};
{\ar "e";"g"};{\ar "hh";"d"};%
{\ar "d";"g"};
{\ar "g";"h"};{\ar "g";"hh"};%
{\ar "h";"e"};{\ar "ee";"hh"};%
{\ar "k";"g"};
{\ar "kk";"hh"};{\ar "hh";"k"};%
{\ar "ff";"hhh"};{\ar "hhh";"jj"};{\ar "hhh";"hh"};
{\ar "hh";"ff"};{\ar "hh";"jj"};
\end{xy}}
\caption{Quiver $Q^{\leq \xi}_2$ in type $A_4$, where $\xi(1)=0$, $\xi(2)=-1$, $\xi(3)=-2$ and $\xi(4)=-1$.} \label{quiver in type A4}
\end{figure}
\end{example}

Recall that for any $l\in \mathbb{Z}_{\geq1}$, the quiver $\Gamma^{\leq \xi}_l$ is defined in Section \ref{quivers and cluster algebras}. We equip $I$ with an order induced by $\xi$:
\begin{align}\label{oder I}
i \geq j \quad \text{if and only if} \quad \xi(i) - \xi(j) \in \{0, 1\}.
\end{align}
Fix an ordering $(i_1, \dots, i_n)$ of $I$ compatible with (\ref{oder I}). For $l \in \mathbb{Z}_{\geq 1}$, define a mutation sequence $\mathscr{S}_l$ in $\Gamma^{\leq \xi}_{l}$ as: $\mathscr{S}_1 = \emptyset$ and for $l\geq 2$,
\begin{align*}
\mathscr{S}_l = & \{ (i_1, \xi(i_1)), (i_1,\xi(i_1)-2), \ldots,  (i_1,\xi(i_1)-2l+4), \\
& \,\,\, (i_2, \xi(i_2)), (i_2,\xi(i_2)-2), \ldots,  (i_2,\xi(i_2)-2l+4), \\
& \qquad \qquad \qquad \qquad \quad \cdots, \cdots, \cdots, \\
& \,\,\, (i_n, \xi(i_n)), (i_n,\xi(i_n)-2), \ldots,  (i_n,\xi(i_n)-2l+4) \}.
\end{align*}

The isomorphism in Theorem \ref{cluster structures on C} induces the following correspondence:
\begin{align*}
\left[W^{(i)}_{l-1,\xi(i)-2l+4}\right] & \longmapsto z_{i,\xi(i)-2l+4}, \\
\left[W^{(i)}_{l,\xi(i)-2l+2}\right]  & \longmapsto  z_{i,\xi(i)-2l+2}, \\
\left[W^{(i)}_{l+1,\xi(i)-2l}\right] & \longmapsto z_{i,\xi(i)-2l}.
\end{align*}

Let $\mathbf{x}'$ be the cluster obtained by applying the mutation sequence $\mathscr{S}_l$ to the initial cluster and denote by $\mu_{\mathscr{S}_l}(x)$ the cluster variable corresponding to the position of $x$ in this cluster. Following \cite[Section 2.2.3 and Theorem 3.1~(b)]{HL16}, each exchange relation in the mutation sequence $\mathscr{S}_l$ corresponds to a T-system equation. Thus we have the following correspondence:
\begin{align*}
\mu_{\mathscr{S}_l}\left(\left[W^{(i)}_{l-1,\xi(i)-2l+4}\right]\right) = \left[W^{(i)}_{l-1,\xi(i)-2l+2}\right]  \longmapsto \mu_{\mathscr{S}_l}(z_{i,\xi(i)-2l+4}).
\end{align*}
It follows from Definition \ref{defintion of the quiver} that the principal quiver of $Q^{\leq \xi}_l$ is naturally isomorphic to $Q_{\xi}$, and $Q^{\leq \xi}_l$ is a subquiver of $\mu_{\mathscr{S}_l}(\Gamma^{\leq \xi}_l)$ with vertices:
\begin{align*}
\left\{ (i, \xi(i) - 2l + 4),  (i, \xi(i) - 2l + 2), (i, \xi(i) - 2l) \mid i \in I \right\}.
\end{align*}

For simplicity, let $\mathcal{A}^{\leq \xi}_l$ be the cluster subalgebra of $K_0(\mathscr{C}^{\leq \xi}_l)$ with initial seed
\begin{align*}
\left( \left([W^{(i)}_{l-1,\xi(i)-2l+2}],  [W^{(i)}_{l,\xi(i)-2l+2}], [W^{(i)}_{l+1,\xi(i)-2l}] \right)_{i\in I}, Q^{\leq \xi}_l\right),
\end{align*}
where $W^{(i)}_{0,r_i} = 1$ when $r_i > \xi(i)$, and $\left\{ [W^{(i)}_{l-1,\xi(i)-2l+2}], [W^{(i)}_{l+1,\xi(i)-2l}] \mid i \in I \right\}$ is the set of coefficients. Thus we obtain the following proposition.

\begin{proposition}\label{prop 4.6}
For any $l \in \mathbb{Z}_{\geq 1}$, the map
\begin{align*}
[W^{(i)}_{l-1,\xi(i)-2l+2}]  & \longmapsto \mu_{\mathscr{S}_l}(z_{i,\xi(i)-2l+4})  \\
[W^{(i)}_{l,\xi(i)-2l+2}]  & \longmapsto  z_{i,\xi(i)-2l+2} \\
[W^{(i)}_{l+1,\xi(i)-2l}] & \longmapsto  z_{i,\xi(i)-2l}
\end{align*}
induces a cluster algebra isomorphism $\mathcal{A}^{\leq \xi}_l \cong \mathcal{A}(Q^{\leq \xi}_l)$.
\end{proposition}

Recall the definition of cluster modules in Definition \ref{a new family of real prime modules}. In the following we present some properties of the cluster modules.

\begin{theorem}\label{indecomposable are real prime}
The map $\Psi^{\leq \xi}_l$ reduces to a map from the set of indecomposable objects in $\mathcal{C}_{Q_{\xi}}$ to the set of real prime modules in $\mathscr{C}^{\leq \xi}_l$.
\end{theorem}

We will prove Theorem \ref{indecomposable are real prime} in Section \ref{proof of Theorem 4.7}.      
As a result of Theorem \ref{indecomposable are real prime}, for $l = 1$, our parameters satisfy
\[
u_i(1) = 1 \text{ for all $i \in I$}, \,\, \widetilde{\textbf{x}}^{\leq \xi}_1 = \textbf{z}^{\leq \xi}_1, \,\, \text{and } \mathcal{C}_{Q_\xi} = \mathcal{C}^{\leq \xi}_1.
\] 
The map $\Psi^{\leq \xi}_1$ induces the bijection
\begin{align*}
\Phi^{\leq \xi}_1: \left\{ \makecell{\text{indecomposable} \\ \text{objects in $\mathcal{C}^{\leq \xi}_1$}} \right\} & \longrightarrow  \left\{ \makecell{ \text{the classes of non-frozen} \\ \text{real prime modules in $\mathscr{C}^{\leq \xi}_1$}} \right\},
\end{align*}
In this case, the images of the indecomposable objects in $\mathcal{C}^{\leq \xi}_1$ under the map $\Psi^{\leq \xi}_1$ recover Hernandez--Leclerc modules \cite{BC19,DS25}. So we have the following corollary.

\begin{corollary}\label{HL are cluster modules}
Every Hernandez--Leclerc module is a prime cluster module.
\end{corollary}

By applying Theorem \ref{indecomposable are real prime}, we show in Theorem \ref{cluster modules are reachable} that cluster modules are reachable and real.

\begin{theorem}\label{cluster modules are reachable}
The cluster modules are reachable real modules.
\end{theorem}
\begin{proof}
Since $Q_{\xi}$ is the principle quiver of $\check{Q}_\xi$, an arbitrary rigid object in $\mathcal{C}_{Q_{\xi}}$ corresponds to a cluster monomial in $\mathcal{A}(\check{Q}_\xi)$. From the proof of Theorem \ref{indecomposable are real prime}, it follows that a cluster monomial in $\mathcal{A}(\check{Q}_\xi)$ corresponds to a cluster monomial in $\mathcal{A}^{\leq \xi}_l\subset K_0(\mathscr{C}^{\leq \xi}_l)$. Hence an arbitrary rigid object in $\mathcal{C}_{Q_{\xi}}$ corresponds to a cluster monomial in $\mathcal{A}^{\leq \xi}_l\subset K_0(\mathscr{C}^{\leq \xi}_l)$ for some $l\in \mathbb{Z}_{\geq 1}$. Therefore, the cluster modules are reachable real objects.
\end{proof}

\begin{corollary}
Let $M_1, M_2$ be two indecomposable objects in $\mathcal{C}_{Q_{\xi}}$. Then $\Psi^{\leq \xi}_{l}(M_1)\otimes \Psi^{\leq \xi}_{l}(M_2)$ is real if and only if $M_1\oplus M_2$ is rigid in $\mathcal{C}_{Q_{\xi}}$.
\end{corollary}

In Theorem \ref{KR are Cluster module}, we show that the cluster modules contain Kirillov--Reshetikhin modules as a subfamily. The proof of Theorem  \ref{KR are Cluster module} appears in Section \ref{proof of KR are Cluster module}.

\begin{theorem}\label{KR are Cluster module}
Every Kirillov--Reshetikhin module is a cluster module.
\end{theorem}

\section{System of equations for cluster modules}\label{a new family of real prime modules and Equational systems}

\subsection{Notation}

Let $Q$ be a Dynkin quiver, i.e., an orientation of a Dynkin diagram, and $\mathcal{C}_Q$ the cluster category associated to $Q$. Following \cite{BMR07}, there is an equivalence of categories:
\begin{align*}
\mathscr{F}_{Q}= \textup{Hom}_{\mathcal{C}_Q}(\mathbb{K}Q,-): \mathcal{C}_Q/\text{add} \mathbb{K}Q [1] \overset{\sim}{\to} \mathbb{K}Q\text{-mod},
\end{align*}
where $\mathbb{K}Q\text{-mod}$ denotes the category of the finitely generated $\mathbb{K}Q$-modules.

Assume that $(L,N)$ is an exchange pair in $\mathcal{C}_Q$ with exchange triangles
\begin{align} \label{non-split triangles}
L\to M \to N\to L[1],  \quad  N \to M' \to L \to N[1].
\end{align}
By a result of \cite[Section 5.2]{DS25}, we always assume without loss of generality that
\[
\mathbf{g}(M)=\mathbf{g}(L)+\mathbf{g}(N)
\]
(exchanging $L$ and $N$ if necessary) in the non-split exchangeable triangles (\ref{non-split triangles}). Then by \cite[Theorem 68 and Lemma 31~(2)]{CIEFR21} (or see \cite[Theorem 1.5]{ST16})
\begin{align}\label{cluster exchanges in a cluster category}
X_L X_{N}=X_M + \textbf{y}^{\cc}X_{M'},
\end{align}
where the vector $\cc$ is defined as the dimension vector of the image of the morphism $h:\tau^{-1}L \to N$. This morphism $h$ is unique up to scaling, because 
\[
\textup{Hom}(\tau^{-1}L,N)\cong D \textup{Ext}^1(N,L)= \mathbb{K},
\]
where the first identity holds by the Auslander--Reiten formula in the cluster category and the second identity is true because $L$ and $N$ form an exchange pair.

For any $M\in \mathbb{K}Q\text{-mod}$, the socle $\text{soc}(M)$  of $M$ is the submodule generated by all simple submodules of $M$. We define
\begin{align}\label{scalar kappa}
\kappa (L,M,N) = \underline{\textup{dim}}(\textup{soc}(\mathscr{F}_{Q}L))+\underline{\textup{dim}}(\textup{soc}(\mathscr{F}_{Q}N))-\underline{\textup{dim}}(\textup{soc}(\mathscr{F}_{Q} M)).
\end{align}
Similarly, we define $\kappa(L,M',N)$.

Following \cite[Theorem 5.11 and Theorem 5.16]{DS25}, the vectors 
\[
\kappa (L,M,N), \quad \kappa (L,M',N)+\mathbf{g}(\textup{Im}(h))
\]
are non-negative integer vectors.

\subsection{Cluster algebra $\mathcal{A}(\check{Q}_\xi)$}\label{Cluster exchange relations I}

Recall that $\xi$ is a height function. There is a Dynkin quiver $Q_\xi$, which is defined by
\[
\text{$i \to j$ if $\xi(i)=\xi(j)+1$}.
\]

Let $\check{Q}_\xi$ be a quiver with $n$ mutable vertices, still denoted by $I=\{1,2,\ldots,n\}$, and $n$ frozen vertices, denoted by $I'=\{1',2',\ldots,n'\}$, and with arrows given by the following rule:
\begin{itemize}
\item[(1)] $i \to j$ if and only if $\xi(i)=\xi(j)+1$, and $C_{ij}\neq 0$;
\item[(2)] For each $i$, there is an arrow $i'\to i$;
\item[(3)] For each arrow $j \to i$, there is an arrow $i\to j'$.
\end{itemize}
By definition, $Q_\xi$ is the principal quiver of $\check{Q}_\xi$.

Let $\mathbb{P}=\text{Trop}(f_1,\ldots,f_n)$ be a tropical semifield, and for any $j\in I$
\begin{align}\label{equation y f}
y_j= f^{-1}_j \prod_{i: i\to j} f_i, \quad  \widehat{y}_j = y_j \left(\prod_{i: i\to j} x^{-1}_i  \prod_{i:j\to i} x_i \right).
\end{align}
Denote by $\mathcal{A}(\check{Q}_\xi)$ the cluster algebra with initial seed $(\textbf{x},\textbf{y},\check{Q}_\xi)$, where $\textbf{x}=(x_i)_{i\in I}$, $\textbf{y}=(y_i)_{i\in I}$, and $\check{Q}_\xi$ is defined as above.

The goal of the section is to give an explicit description of the exchange relations in the cluster algebra
$\mathcal{A}(\check{Q}_\xi)$.  We start by recalling a few results that we shall need.

For vectors ${\bf f}=(f_i)_{i\in I}$ and $d\in \mathbb{Z}^I$, let ${\bf f}^{d}=\prod_{i\in I} f^{d_i}_i$.

\begin{lemma}[{\cite[Lemma 5.9 and Equation (5.15)]{DS25}}] \label{tropical F-polynomials quiver 1}
Let $M \in \mathcal{C}_{Q_\xi}$. Then
\[
F_M |_{\mathbb{P}} (y_1,\ldots,y_{n}) = {\bf f}^{-\underline{\textup{dim}}(\textup{soc}(\mathscr{F}_{Q_\xi} M))}.
\]
\end{lemma}

Let $N \mapsto \check X_N$ be the bijection from the set of indecomposable objects in $\mathcal{C}_{Q_{\xi}}$ to the set of cluster variables in $\mathcal{A}(\check{Q}_\xi)$, such that
\begin{align*} 
\check X_N = \frac{\textbf{x}^{\mathbf{g}(N)} }{F_{N}\vert_{\mathbb{P}}(y_1,\ldots,y_n)} F_{N}(\hat{y}_1,\ldots,\hat{y}_n),
\end{align*}
where $y_i$, $\hat y_i$ are defined in Equation (\ref{equation y f}).

\begin{proposition}\label{an expression of extended g-vectors}
Let $\widetilde{\textbf{x}}=(\textbf{x}\,,\textbf{f}\,)$, and $N$ an indecomposable object in $\mathcal{C}_{Q_{\xi}}$. Then
\begin{align*}
\widetilde{\textbf{x}\,}^{\widetilde{\mathbf{g}}(N)} = \frac{\textbf{x}\,^{\mathbf{g}(N)}}{F_{N}\vert_{\mathbb{P}}(y_1,\ldots,y_n)}  = \textbf{x}\,^{\mathbf{g}(N)} {\bf f\,}^{\underline{\textup{dim}}(\textup{soc}(\mathscr{F}_{Q_\xi}N))}.
\end{align*}
\end{proposition}
\begin{proof}
It follows from Lemma \ref{tropical F-polynomials quiver 1}.
\end{proof}

As a consequence, Proposition \ref{an expression of extended g-vectors} gives an algorithm for the extended $\mathbf{g}$-vectors in terms of quiver representation.

By Equation (\ref{cluster exchanges in a cluster category}), if we specialize all initial cluster variables to 1, the following identity holds:
\begin{align} \label{F-polynomials identity}
F_L(y_1,\ldots,y_n)  F_N(y_1,\ldots,y_n) =  F_M(y_1,\ldots,y_n) + \textbf{y}^{\underline{\text{dim}}(\textup{Im}(h))} F_{M'}(y_1,\ldots,y_n).
\end{align}

The following lemma is a reformulation of Lemma 5.13 of \cite{DS25}. For the convenience of the reader we include a proof.

\begin{lemma} \label{interpretation of haty}
Let $M$ be a $\mathbb{K}Q_\xi$-module. Then
\begin{align*}
\hat{\textbf{y}}^{\underline{\dim}(M)} = \prod_{i \in I} \hat{y}_i^{\underline{\dim}(M_i)} =  {\bf x}^{a(M)}  {\bf f}^{\mathbf{g}(M)},
\end{align*}
where $a(M)\in \mathbb{Z}^I$ is defined by
\begin{align*}
a(M)_i=\sum_{j:j\to i} \underline{\textup{dim}}(M_j) - \sum_{j:i\to j} \underline{\textup{dim}}(M_j).
\end{align*}
\end{lemma}
\begin{proof}
It follows from Lemma 2.3 in \cite{Pal08} that the $j$-component of the $\mathbf{g}$-vector $\mathbf{g}(M)$ of $M$ is
\[
\underline{\text{dim}}(M_j)-\sum_{i:j\to i}\underline{\text{dim}}(M_i).
\]
Hence
\begin{align*}
\hat{\textbf{y}}^{\underline{\dim}(M)} & = \prod_{j \in I} \left( f^{-1}_j \prod_{i: i\to j} f_i x^{-1}_i  \prod_{i:j\to i} x_i \right)^{\underline{\dim}(M_j)} = {\bf x}^{a(M)}  {\bf f}^{\mathbf{g}(M)}.
\end{align*}
\end{proof}

We now present the main result of this section. It gives an explicit description of all exchange relations in the cluster algebra $\mathcal{A}(\check{Q}_\xi)$.

\begin{theorem}\label{all cluster relations associated to quiver I}
Let $(L,N)$ be an exchange pair in the cluster category $\mathcal{C}_{Q_{\xi}}$ with exchange triangles (\ref{non-split triangles}). The corresponding exchange relation in $\mathcal{A}(\check{Q}_\xi)$ has the following form
\begin{align*}
\check X_L  \check X_N  =  \check X_M  {\bf f}^{\kappa (L,M,N)} +  \check X_{M'}  {\bf f}^{\kappa (L,M',N)+\mathbf{g}(\textup{Im}(h))},
\end{align*}
where $\kappa (L,M,N)$ and  $\kappa (L,M',N)$ are defined as Equation (\ref{scalar kappa}). The coefficient free version of this formula is well-known from \cite{BMRRT06}.
\end{theorem}

\begin{proof}
We assume without loss of generality that $\mathbf{g}(M)=\mathbf{g}(L)+\mathbf{g}(N)$. By Theorem \ref{Fomin-Zelevinsky separation}, for any indecomposable module $N$, we have
\begin{align}\label{separation formula of XN}
\check X_N  = \widetilde{\textbf{x}}^{\widetilde{\mathbf{g}}(N)} F_N (\hat{y}_1,\ldots,\hat{y}_n),
\end{align}
and $\mathbf{g}(M') = \mathbf{g}(M)+a(\textup{Im}(h))$ by \cite[Lemma 5.14]{DS25}.

By Equation (\ref{F-polynomials identity}), Proposition \ref{an expression of extended g-vectors} and Lemma \ref{interpretation of haty}, we have
\begin{align*}
\check X_L  \check X_N & = \widetilde{\textbf{x}}^{{\widetilde{\mathbf{g}}(L)}+{\widetilde{\mathbf{g}}(N)}} F_L(\hat{y}_1,\ldots,\hat{y}_n) F_N(\hat{y}_1,\ldots,\hat{y}_n) \\
& = \textbf{x}^{\mathbf{g}(L)+\mathbf{g}(N)} {\bf f}^{\underline{\textup{dim}}(\textup{soc}(\mathscr{F}_{Q_\xi} L))+\underline{\textup{dim}}(\textup{soc}(\mathscr{F}_{Q_\xi} N))}  \\
& \quad  \times \left( F_{M}(\hat{y}_1,\ldots,\hat{y}_n) +\hat{\textbf{y}}^{\underline{\text{dim}}(\textup{Im}(h))}F_{M'}(\hat{y}_1,\ldots,\hat{y}_n) \right) \\
& =  \textbf{x}^{\mathbf{g}(M)} F_{M}(\hat{y}_1,\ldots,\hat{y}_n) {\bf f}^{\underline{\textup{dim}}(\textup{soc}(\mathscr{F}_{Q_{\xi}} L))+\underline{\textup{dim}}(\textup{soc}(\mathscr{F}_{Q_{\xi}} N))} \\
& \quad  + \textbf{x}^{\mathbf{g}(M)} {\bf f}^{\underline{\textup{dim}}(\textup{soc}(\mathscr{F}_{Q_\xi} L))+\underline{\textup{dim}}(\textup{soc}(\mathscr{F}_{Q_{\xi}} N))}  \textbf{x}^{a(\textup{Im}(h))} {\bf f}^{\mathbf{g}(\textup{Im}(h))} F_{M'}(\hat{y}_1,\ldots,\hat{y}_n) \\
&  = \check X_M {\bf f}^{\kappa (L,M,N)} + \textbf{x}^{\mathbf{g}(M')} F_{M'}(\hat{y}_1,\ldots,\hat{y}_n)  {\bf f}^{\underline{\textup{dim}}(\textup{soc}(\mathscr{F}_{Q_{\xi}} L)+\underline{\textup{dim}}(\textup{soc}(\mathscr{F}_{Q_{\xi}} N))) + \mathbf{g}(\textup{Im}(h))}\\
&  = \check X_M {\bf f}^{\kappa (L,M,N)} + \check X_{M'} {\bf f}^{\kappa (L,M',N)+\mathbf{g}(\textup{Im}(h))},
\end{align*}
where the second equation uses Proposition \ref{an expression of extended g-vectors} and Equation (\ref{F-polynomials identity}), the third equation uses Lemma \ref{interpretation of haty}, and the last two equations follow from (\ref{scalar kappa}) and (\ref{separation formula of XN}). This completes the proof.
\end{proof}

We present an example that illustrates Theorem \ref{all cluster relations associated to quiver I}.
\begin{example}
Let $L=S(3)$ and $N=\substack{\,1\, \\ 2}$ from the previous Example \ref{computation of the extended g-vectors}.  Then by definition
\begin{align*}
& M = \substack{\,1 \, \\ \,2\, \\ 3}, \quad M'= \substack{1}, \\
& \kappa (L,M,N) =  \left(\begin{array}{c} 0 \\ 1 \\ 0 \end{array} \right), \\
& \kappa (L,M',N)+\mathbf{g}(\textup{Im}(h)) = \left(\begin{array}{c} -1 \\ 1 \\ 1 \end{array} \right) +  \left(\begin{array}{c} 1 \\ -1 \\ 0 \end{array} \right) = \left(\begin{array}{c} 0 \\ 0 \\ 1 \end{array} \right).
\end{align*}
By Theorem \ref{all cluster relations associated to quiver I}, we have
\[
\check X_L  \check X_N  =  \check X_M  f_2  +  \check X_{M'} f_3.
\]
\end{example}

\subsection{System of equations for cluster modules}

In this section, we present a system of equations satisfied by $q$-characters of cluster modules, which naturally generalizes the relations for Hernandez--Leclerc modules in \cite[Theorem 5.8]{DS25}. 

\begin{theorem}\label{cluster exchange relations}
Let $(L,N)$ be an exchange pair in $\mathcal{C}_{Q_{\xi}}$ with exchange triangles (\ref{non-split triangles}). Then the following equation holds in $K_0(\mathscr{C}^{\leq \xi}_l)$:
\begin{align}\label{cluster exchange I in the subcategory}
\begin{split}
\left[ \Psi^{\leq \xi}_{l}(L) \right] \left[ \Psi^{\leq \xi}_{l}(N) \right] & = \left[ \Psi^{\leq \xi}_{l}(M) \right] \left(\prod_{i \in I} \left[ L(u_i(l)) \right]^{\alpha_i} \left[ L(v_i(l)) \right]^{\alpha_i} \right) \\
& \quad + \left[ \Psi^{\leq \xi}_{l}(M') \right] \left( \prod_{j \in I} \left[ L(u_j(l)) \right]^{\beta_j} \left[ L(v_j(l)) \right]^{\beta_j} \right),
\end{split}
\end{align}
where $l\in \mathbb{Z}_{\geq 1}$, the vectors $(\alpha_i)_{i\in I}$ and $(\beta_i)_{i\in I}$ are given by
\begin{align*}
(\alpha_i)_{i\in I} & =  \kappa (L,M,N), \\
(\beta_j)_{j\in I} & = \kappa (L,M',N) + \mathbf{g}(\operatorname{Im}(h)).
\end{align*}
Moreover,
\[
\hw \left(\Psi^{\leq \xi}_{l}(L) \otimes \Psi^{\leq \xi}_{l}(N)\right) = \hw \left( \Psi^{\leq \xi}_{l}(M) \otimes \left(\bigotimes_{i \in I} \left(L(u_i(l)) \otimes L(v_i(l)) \right)^{\alpha_i} \right)\right).
\]
\end{theorem}

\begin{proof}
By Proposition \ref{prop 4.6}, we identify the coefficients and the cluster variables
\begin{align*}
[L(u_i(l))] [L(v_i(l))] & \longmapsto  f_i, \\
[W^{(i)}_{l,\xi(i)-2l+2}] & \longmapsto  x_i,
\end{align*}
then the cluster structure of $\mathcal{A}^{\leq \xi}_l$ is the same as the cluster structure of $\mathcal{A}(\check{Q}_\xi)$. The formula (\ref{cluster exchange I in the subcategory}) follows from Theorem \ref{all cluster relations associated to quiver I} via the substitution 
\[
f_i \longmapsto [L(u_i(l))] [L(v_i(l))] \text{ for $i \in I$}, \quad \check X_L  \longmapsto \left[ \Psi^{\leq \xi}_{l}(L) \right].
\]

%

By the definition of $\Psi^{\leq \xi}_l$, we have
\begin{align*}
\hw \left(\Psi^{\leq \xi}_{l}(L) \right) &= (\textbf{x}^{\leq \xi}_l)^{\mathbf{g}(L)} ({\bf c}^{\leq \xi}_l)^{\underline{\dim}(\operatorname{soc}(\mathscr{F}_{Q_\xi} L))}, \\
\hw \left(\Psi^{\leq \xi}_{l}(N) \right) &= (\textbf{x}^{\leq \xi}_l )^{\mathbf{g}(N)} ({\bf c}^{\leq \xi}_l)^{\underline{\dim}(\operatorname{soc}(\mathscr{F}_{Q_\xi} N))},
\end{align*}
and
\begin{align*}
& \hw \left(\Psi^{\leq \xi}_{l}(M) \otimes \bigotimes_{i \in I} \left(L(u_i(l)) \otimes L(v_i(l)) \right)^{\otimes \alpha_i} \right) \\
&= \hw \left(\Psi^{\leq \xi}_{l}(M)\right) \hw \left( \bigotimes_{i \in I} \left(L(u_i(l)) \otimes L(v_i(l)) \right)^{\otimes \alpha_i} \right) \\
&= (\textbf{x}^{\leq \xi}_l)^{\mathbf{g}(M)} ({\bf c}^{\leq \xi}_l)^{\underline{\dim}(\operatorname{soc}(\mathscr{F}_{Q_\xi} M)) + \kappa (L,M,N)}.
\end{align*}
Therefore, by the definition of $\kappa (L,M,N)$,
\[
\hw \left(\Psi^{\leq \xi}_{l}(L) \otimes \Psi^{\leq \xi}_{l}(N)\right) = \hw \left( \Psi^{\leq \xi}_{l}(M) \otimes \bigotimes_{i \in I} \left(L(u_i(l)) \otimes L(v_i(l)) \right)^{\otimes \alpha_i} \right).\]
\end{proof}

When $l=1$, Equation (\ref{cluster exchange I in the subcategory}) in Theorem \ref{cluster exchange relations} reduces to the formula established in \cite[Theorem 5.8]{DS25}. Moreover, Equation (\ref{cluster exchange I in the subcategory}) implies the existence of a pair of exact sequences in $\mathscr{C}^{\leq \xi}_l$:

\begin{align*}
0 \to \Psi^{\leq \xi}_{l}(M') \otimes \bigotimes_{j \in I} (L(u_j(l)) & \otimes  L(v_j(l)))^{\otimes \beta_j}  \to  \Psi^{\leq \xi}_{l}(L) \otimes  \Psi^{\leq \xi}_{l}(N)  \\
& \to  \Psi^{\leq \xi}_{l}(M) \otimes \bigotimes_{i \in I} \left(L(u_i(l)) \otimes L(v_i(l)) \right)^{\otimes \alpha_i}  \to 0, \\
0 \to  \Psi^{\leq \xi}_{l}(M) \otimes \bigotimes_{i \in I}  (L(u_i(l)) & \otimes  L(v_i(l)))^{\otimes \alpha_i}  \to  \Psi^{\leq \xi}_{l}(N) \otimes  \Psi^{\leq \xi}_{l}(L)  \\
& \to \Psi^{\leq \xi}_{l}(M') \otimes \bigotimes_{j \in I}  \left(L(u_j(l)) \otimes L(v_j(l)) \right)^{\otimes \beta_j} \to 0,
\end{align*}
or the same sequences with arrows reversed.

\section{Proofs of Theorem \ref{indecomposable are real prime} and Theorem \ref{KR are Cluster module}}\label{Proofs of two theorems}

In this section, we prove Theorem \ref{indecomposable are real prime} and Theorem \ref{KR are Cluster module}.

\subsection{Proof of Theorem \ref{indecomposable are real prime}}\label{proof of Theorem 4.7}
It suffices to prove that $\Psi^{\leq \xi}_{l}(M)$ is a real prime module for an arbitrary indecomposable object $M\in \mathcal{C}_{Q_{\xi}}$.

Note that the vertices $i'$ (if present) and $i''$ in $Q^{\leq \xi}_l$ are symmetric with respect to $Q_\xi$. Identifying each pair of vertices $i'$ and $i''$ (when they exist) as a single vertex, that is, for $i \in I$,
\begin{align*}
\mu_{\mathscr{S}_l}(z_{i,\xi(i)-2l+4}) \times z_{i,\xi(i)-2l} &  \longmapsto  f_i, \\
z_{i,\xi(i)-2l+2} &  \longmapsto  x_i,
\end{align*}
the cluster structure of $\mathcal{A}(Q^{\leq \xi}_l)$ is the same as the cluster structure of $\mathcal{A}(\check{Q}_\xi)$. By Proposition \ref{prop 4.6}, we identify the coefficients and the cluster variables
\begin{align*}
[
W^{(i)}_{l-1,\xi(i)-2l+2}] [W^{(i)}_{l+1,\xi(i)-2l}] & \longmapsto f_i, \\
[W^{(i)}_{l,\xi(i)-2l+2}] & \longmapsto  x_i,
\end{align*}
then the cluster structure of $\mathcal{A}^{\leq \xi}_l$ is the same as the cluster structure of $\mathcal{A}(\check{Q}_\xi)$.

On the other hand, we have
\begin{align*} 
\hw (W^{(i)}_{l-1,\xi(i)-2l+2}) & =  u_i(l),\\
\hw (W^{(i)}_{l,\xi(i)-2l+2}) & = z_{i,\xi(i)-2l+2}, \\
\hw (W^{(i)}_{l+1,\xi(i)-2l}) & = v_i(l).
\end{align*}
Recall that in Section \ref{definition of cluster modules}, 
\begin{align*}
\textbf{x}^{\leq \xi}_l =(z_{i,\xi(i)-2l+2})_{i\in I}, \quad {\bf c}^{\leq \xi}_l = (u_i(l) v_i(l))_{i\in I}, \quad  \widetilde{\textbf{x}}^{\leq \xi}_l = (\textbf{x}^{\leq \xi}_l, {\bf c}^{\leq \xi}_l),
\end{align*}
and $\mathbb{P}_l$ is the tropical semifield generated by $u_i(l)$ and  $v_i(l)$, with $i\in I$. By Proposition \ref{an expression of extended g-vectors}, the extended $\mathbf{g}$-vector of the cluster variable associated to $M$ in $\mathcal{A}(\check{Q}_\xi)$ is
\[
\widetilde{\mathbf{g}}(M) =
\begin{pmatrix}
\mathbf{g}(M) \\
\underline{\dim}(\operatorname{soc}(\mathscr{F}_{Q_{\xi}} M))
\end{pmatrix}.
\]
After identifying the coefficients and the cluster variables
\begin{align*}
u_i(l) v_i(l) \longmapsto f_i,  \quad  z_{i,\xi(i)-2l+2} \longmapsto  x_i,
\end{align*}
the cluster variable in $\mathcal{A}^{\leq \xi}_l$ associated to $M$ is 
\begin{align*} 
(\widetilde{\textbf{x}}^{\leq \xi}_l)^{\widetilde{\mathbf{g}}(M)} F_M \left(\widehat{u_1(l)v_1(l)},\ldots,\widehat{u_n(l)v_n(l)}\right),
\end{align*}
where for any $j \in I$
\[
\widehat{u_j(l)v_j(l)} = \left( \frac{\prod_{i: i\to j \text{ in $Q^{\leq \xi}_l$}} u_i(l)v_i(l) }{(u_j(l)v_j(l))} \right) \left( \frac{\prod_{i:j\to i \text{ in $Q^{\leq \xi}_l$}} z_{i,\xi(i)-2l+2}}{\prod_{i: i\to j \text{ in $Q^{\leq \xi}_l$}} z_{i,\xi(i)-2l+2}} \right).
\]

By the same argument as in the proof of Theorem \ref{the classification of real modules under reachability}, the highest $\ell$-weight monomial of the cluster variable in $\mathcal{A}^{\leq \xi}_l$ associated to $M$ is $(\widetilde{\textbf{x}}^{\leq \xi}_l)^{\widetilde{\mathbf{g}}(M)}$.  Equivalently, 
\[
\left[L\left((\widetilde{\textbf{x}}^{\leq \xi}_l)^{\widetilde{\mathbf{g}}(M)}\right)\right] = [\Psi^{\leq \xi}_{l}(M)]
\]
is a cluster variable in $\mathcal{A}^{\leq \xi}_l$. In \cite{Qin17}, Qin proved that the modules corresponding to cluster monomials are real. Hence $\Psi^{\leq \xi}_{l}(M)$ is real. It was proved in \cite[Corollary 8.6]{GLS13} that the module in $\mathscr{C}^{\leq \xi}_l$  corresponding to a cluster variable is prime. Since $\mathcal{A}^{\leq \xi}_l$ is a cluster subalgebra of $K_0(\mathscr{C}^{\leq \xi}_l)$ and $[\Psi^{\leq \xi}_{l}(M)]$ corresponds to a cluster variable, $\Psi^{\leq \xi}_{l}(M)$ is prime.

Therefore, $\Psi^{\leq \xi}_{l}(M)$ is a real prime module for an arbitrary indecomposable object $M\in \mathcal{C}_{Q_{\xi}}$.

\subsection{Proof of Theorem \ref{KR are Cluster module}}\label{proof of KR are Cluster module}

Since $\mathcal{C}_{Q_{\xi}}$ is the cluster category of the hereditary path algebra of $Q_{\xi}$, for every vertex $i$ there exist the following pairs of non-split exchange triangles:
\begin{align*}
& P(i)[1] \to  0  \to  I(i) \to P(i)[2]=I(i), \\
& I(i) \to \left(\bigoplus_{j:i\to j} P(j)[1]\right) \oplus \left(\bigoplus_{j:j\to i} I(j)\right) \to P(i)[1].
\end{align*}
By Theorem \ref{cluster exchange relations}, we have
\begin{align*}
(\alpha_i)_{i\in I} = \underline{\dim}(S(i)), \quad  (\beta_j)_{j\in I} = 0,
\end{align*}
and
\begin{align}\label{$T$-system equation in our setting}
\begin{split}
\left[ \Psi^{\leq \xi}_{l}(P(i)[1])\right] \left[ \Psi^{\leq \xi}_{l}(I(i))\right] & = \left[ L(u_i(l)) \right] \left[ L(v_i(l)) \right] \\
& + \prod_{j:i\to j} \left[ \Psi^{\leq \xi}_{l}(P(j)[1])\right] \prod_{j:j\to i} \left[\Psi^{\leq \xi}_{l}(I(j))\right].
\end{split}
\end{align}

By Theorem \ref{property of our map}, Equation (\ref{$T$-system equation in our setting}) becomes
\begin{align}\label{5.3}
\begin{split}
\left[ W^{(i)}_{l,\xi(i)-2l+2} \right] \left[ W^{(i)}_{l,\xi(i)-2l}  \right] & = \left[ W^{(i)}_{l-1,\xi(i)-2l+2}  \right]  \left[ W^{(i)}_{l+1,\xi(i)-2l}  \right]  \\
& + \prod_{j:i\to j} \left[ W^{(j)}_{l,\xi(j)-2l+2} \right] \prod_{j:j\to i} \left[ W^{(j)}_{l,\xi(j)-2l} \right].
\end{split}
\end{align}
Following our convention, $i \to j$ if and only if $C_{ij} = -1$ and $\xi(i)=\xi(j)+1$. Thus,
$$ \prod_{j:i\to j} \left[ W^{(j)}_{l,\xi(j)-2l+2} \right]  \prod_{j:j\to i} \left[ W^{(j)}_{l,\xi(j)-2l} \right] = \prod_{j: C_{ij}=-1} \left[ W^{(j)}_{l,\xi(i)-2l+1} \right],$$
and the formula (\ref{5.3}) reduces to an equation from T-systems:
\begin{align*}
\left[ W^{(i)}_{l,\xi(i)-2l+2}\right] \left[ W^{(i)}_{l,\xi(i)-2l} \right] = \left[ W^{(i)}_{l-1,\xi(i)-2l+2} \right]  \left[ W^{(i)}_{l+1,\xi(i)-2l} \right] + \prod_{j: C_{ij}=-1} \left[ W^{(j)}_{l,\xi(i)-2l+1}\right],
\end{align*}
see Equation (\ref{equ:T-system equation}).

Therefore, the system of equations in Theorem \ref{cluster exchange relations} encompasses the classical T-system relations for Kirillov--Reshetikhin modules, and the modules arising from this system of equations include Kirillov--Reshetikhin modules as a special subfamily.


\section*{Acknowledgements}
We would like to express our gratitude to Fan Qin and Yilin Wu for insightful discussions.

\end{sloppypar}
\end{document}